\documentstyle[secteqno]{article}
\author{W.~D.~Evans, D.~J.~Harris, J.~Lang }
\title{The approximation numbers of
Hardy--type operators on trees.}

\newcommand{\rank}{\mathop{\rm rank}\nolimits}

\newcommand{\sgn}{\mathop{\rm sgn}\nolimits}

\newtheorem{theorem}{Theorem}[section]
\newtheorem{lemma}[theorem]{Lemma}
\newtheorem{cor}[theorem]{Corollary}
\newtheorem{defn}[theorem]{Definition}
\newtheorem{prop}[theorem]{Proposition}
\newtheorem{remark}[theorem]{Remark}
\newcommand{\bpr}{\begin{prop}}
\newcommand{\beq}{\begin{equation}}
\newcommand{\bd}{\begin{defn}}
\newcommand{\bex}{\begin{example}}
\newcommand{\bc}{\begin{cor}}
\newcommand{\bl}{\begin{lemma}}
\newcommand{\bt}{\begin{theorem}}
\newcommand{\br}{\begin{remark}}
\newcommand{\epr}{\end{prop}}
\newcommand{\eeq}{\end{equation}}
\newcommand{\ed}{\end{defn}}
\newcommand{\eex}{\end{example}}
\newcommand{\ec}{\end{cor}}
\newcommand{\el}{\end{lemma}}
\newcommand{\et}{\end{theorem}}
\newcommand{\er}{\end{remark}}
\newcommand{\R}{ {\bf R} }
\newcommand{\N}{ {\bf N} }

\newcommand{\K}{ {\bf K} }
\newcommand{\Z}{ {\bf Z} }

\newcommand{\I}{ {\bf I} }
\newcommand{\G}{ {\Gamma} }
\newcommand{\eps}{ {\varepsilon} }
\newcommand{\cP}{ {\cal P} }
\newcommand{\cS}{ {\cal S} }
\newcommand{\cL}{ {\cal L} }
\newcommand{\Ge}{ {\cal G} }
\newcommand{\cF}{ {\cal F} }
\newcommand{\cC}{ {\bf C} }
\newcommand{\Alpha}{ {\cal L} }

\begin{document}
\maketitle

\begin{abstract}
The Hardy operator $T_a$ on a tree $\G$ is defined by
\[(T_af)(x):=v(x) \int^x_a u(t)f(t) dt
\qquad \mbox{ for } a, x\in \G.
\]
Properties of $T_a$ as a map from $L^p(\G)$ into itself are established for
$1\le p \le \infty$. The main result is that, with appropriate assumptions
on $u$ and $v$, the approximation numbers $a_n(T_a)$ of $T_a$ satisfy
\[
(*) \lim_{n\rightarrow \infty} na_n(T_a) = \alpha_p\int_{\G}
|uv|dt
\]
for a specified constant $\alpha_p$ and $1<p<\infty$. This extends
results of Naimark, Newman and Solomyak for $p=2$. Hitherto, for
$p\neq 2$, (*) was unknown even when $\G$ is an interval. Also, upper
and lower estimates for the $l^q$ and weak-$l^q$ norms of $\{a_n(T_a)\}$
are determined.

\end{abstract}

\section{Introduction.}

In [1],[2]and [6] results were established for the Hardy operator
\begin{equation}
Tf(x)=v(x) \int_0^x f(t) u(t) dt
\end{equation}
as a map from $L^p(0,\infty)$ to $L^p(0,\infty)$, for $1\leq p\leq \infty$.
When $p\in(1,\infty)$, it was proved in [2] that under appropriate
conditions on $u$ and $v$ the approximation numbers $a_n(T) $ of
$T$ satisfy
\begin{equation}
\lim_{n\to \infty} n a_n(T) = {1\over \pi} \int_0^{\infty} |u(t)
v(t) | dt
\end{equation}
\footins{1991 Mathematics Subject Classification : 47G10, 47B10}

\noindent when p=2, and when $1<p<\infty$
\begin{eqnarray}
 {\alpha_p \over 4} \int_0^{\infty} |u(t) v(t)| dt & \le &
\liminf_{n\to \infty} n a_n(T)\nonumber \\
& \le & \limsup_{n\to \infty} n
a_n(T) \le  {\alpha_p } \int_0^{\infty} |u(t) v(t)| dt
\end{eqnarray}
for a specified constant $\alpha_p$ depending on $p$. For the cases
$p=\infty$ and $p=1$  similar estimates were derived in [6] but
with $v_s(t)=\lim_{\varepsilon \to 0}
\|v\|_{\infty,(t-\varepsilon, t+\varepsilon)}$ instead of $v(t)$ when
$p=\infty $ and $u_s(t)$ instead of $u(t)$ in the case $p=1$. Furthermore, in [2] and
[6] two-sided estimates are given for the $l^{\alpha}$ and
weak-$l^{\alpha}$ norms of the sequence of approximation numbers
in the case when the Hardy operator is compact.

A special case of the main result in this paper is that the counterpart of
(1.2) in $L^p(0,\infty)$, namely
\begin{equation}
\lim_{n\rightarrow \infty} na_n(T)
=\alpha_p\int_0^{\infty}|u(t)v(t)|dt,
\end{equation}
holds for {\em all } $p\in (1,\infty)$; the general result is that the
analogue of (1.4) when the interval $(0,\infty)$ is replaced by a
tree $\G$ is true.  Such Hardy
operators on trees have already been investigated in [5] where it was shown
that they occur naturally in spectral problems defined on domains
with irregular boundaries. Necessary and sufficient criteria for
the boundedness of Hardy operators between various Lebesgue spaces
on $\G$ are established in [7], but the complex nature of the
problem is such that the neat abstract result is difficult to
apply even for the most elementary of trees. It is therefore to be
expected that the problems of compactness and estimating the
approximation numbers are likely to be much more complicated than
in the interval case. This, confirmed in this paper, is what makes
it so surprising that the analogue of (1.4) for a tree is
established here when $p\neq2$ before it was known for an
interval. Estimates for $l^q$ and weak-$l^q$ norms of the
approximation numbers of $T$ are also obtained.

In [9] the case $p=2$ of the problem subsequently studied in [2]
and [6] for general $p\in[1,\infty]$ was considered and (1.2) proved, using Hilbert
space methods which do not extend to general values of $p$.
The same problem on a tree $\G$ is the subject of [8] where an intensive
study is made of problems on
trees which are closely related to those here, but in the case
$p=2$ only, and using methods which are very different from those in this
paper. The conditions imposed to ensure the validity of the analogue of (1.4)
for a tree in [8] are similar to those here, but a comparison seems difficult
in general (see Remark 6.12). The main difference is that in [8] they relate
to arbitrary partitions of $\G$ into intervals, whereas our partitions are
into connected subsets specifically determined by functions which have a
fundamental role in the analysis.

\section{Preliminaries.}

In this section we recall the definition of a tree $\G$, introduce
a Hardy--type operator on the tree and quote from [7] the
criterion for the boundedness of the operator as a map from
$L^p(\G)$ into $L^p(\G)$.

A tree $\G$ is a connected graph without loops or cycles, where
the edges are non-degenerate closed line segments whose end--points
are the vertices. Each vertex of $\G$ is of finite degree, i.e. only a
finite number of edges emanate from each vertex. For every $x, y
\in \G$ there is a unique polygonal path in $\G$ which joins $x$
and $y$. The distance between $x$ and $y$ is defined to be the
length of this polygonal path and in this way $\G$ is endowed with
a metric topology.

\bl
Let $\tau (\G)$ be the metric topology on $\G$. Then
\begin{description}
\item[(i)] a set $A\subset \G$ is compact if and only if it is closed
and meets only a finite number of edges;
\item[(ii)] $\tau (\G)$ is locally compact;
\item[(iii)] $\G$ is the union of a countable number of edges. Thus if $\G$ is endowed with the natural 1-dimensional Lebesgue
measure it is a $\sigma$-finite measure space.
\end{description}
\el
{\bf Proof.} See [5].
{$\Box $}

Let $x,y \in \G$ and denote by $(x,y)$ the unique path joining $x,y$
in $\G$. For $a\in \G$ we define $t \succeq_a x$ (or $x
\preceq_a t$) to mean that $x$ lies on the path $(a,t)$.
We write $x \prec_a t$ for $x \preceq_a t$ and $x \neq t$. This is
a partial ordering on $\G$ and the ordered graph so formed is
referred to us a tree rooted at $a$ and denoted by $\G_a$. If $a$
is not a vertex we make it one by replacing the edge on which it
lies by two edges. In this way $\G_a$ is the unique finite union
of subtrees $\G_{a,i}$ which intersect only at $a$.

Note  that if $x\notin (a,b)$ then $x \preceq_a y$ if and only if
$x \preceq_b y$.

We shall use the following notation. For a subtree $K$ of $\G$,
$V(K)$, $E(K)$ will denote respectively the sets of vertices,
edges of $K$ and $\partial K$ will denote the set of
boundary points of $K$ in $\G$. The notation $K \subset \subset \G $
will be used to mean
that the closure of $K$ is a compact subset of $\G$; note that,
from Lemma 2.1 (i) this implies that $K$ meets only a finite
number of edges of $\G$. The characteristic function of a set $E$
will be denoted by $\chi_E$. The integral is interpreted in the
following sense :
\[
\int_E g = \sum_e \int_{e\cap E} g
\]
where
\[
\int_{e\cap E} g = \int_c^d g(x)\chi_E(x)dx,
\]
the integral $\int_c^d$ being over the set of points lying in the
path $(c,d)$. For a measurable subset $K$ of $\G$ we define the norm
\[
\|f\|_{p,K} = \big( \int_K |f|^p\big)^{1/p}
\]
on $L^p(K)$. The $L^p(\G)$ norm will be denoted by $\|\cdot\|_p$ if there
is no chance of confusion. Also, if the value of $p$ is clear from the
context, we shall write $\|\cdot\|_K, \|\cdot\|$ for the $
L^p$ norms on
$K,\G$ respectively.  If $A$ is a bounded map between normed spaces $X,Y$
we denote its norm by
\[
\|A|X \rightarrow Y\|
\]
This will be simplified to $\|A\|$ if the spaces $X,Y$ are unambiguous.

A connected subset of $\G$ is a subtree if we add its boundary
points to the set of vertices of $\G$, and hence form new edges
from existing ones. Hereafter we shall always adopt this
convention when we refer to subtrees.

\bd
Let $K$ be a subtree of $\G$ containing $a$. A point $t\in
\partial K$ is said to be maximal if every $x \succ_a t$ lies in
$\G \setminus K$. We denote by $\I_a(\G)$ (or simply $\I_a$) the
set of subtrees $K$ of $\G$ containing a whose boundary points
are all maximal.
\ed

We assume throughout, unless mentioned otherwise, that $u,v$
satisfy the following conditions:
\begin{equation} \label{2.1}
u\in L^{p'}(K), \  v\in L^p(\G ), \ \mbox{ for every
} K\subset \subset \G.
\end{equation}

We may assume, without loss of generality, that $u, v \ge 0$. This
is because multiplication by $\sgn u$ and $\sgn v$ are isometries on
$L^p(\G)$; recall that $\sgn u =u/|u|$ when $u\neq 0$ and $1$
otherwise.

\bd Let $\G$ be a tree, $1\le p \le \infty$, and let $u$ and $v$
be measurable functions on $\G$ which satisfy (\ref{2.1}). For
$x\in \G$ and $f\in L^p(\G)$ we define the Hardy operator by
\begin{equation} \label{2.2}
T_af(x):=v(x)\int_a^x f(t)u(t) dt,\;\;\;a\in \G.
\end{equation}
\ed

In [7] the following necessary and sufficient condition for the
boundedness of $T_a$ was obtained.

\bt
Let $1\le p \le \infty$, $a\in \G$, and suppose $u$ and $v$ satisfy
(\ref{2.1}). For $K\in \I_a$ define
\begin{equation} \label{2.3}
\alpha_K:= \inf\{\|f\|_p: \int_a^t |f| |u| =1 \  \mbox{ for all
} \  t\in \partial K \}.
\end{equation}

Then $T_a$ is bounded from $L^p(\G)$ into $L^p(\G)$ if and only if

\begin{equation} \label{2.4}
A:= \sup_{K\in \I_a} {\|v \chi_{\G \setminus K} \|_p \over
\alpha_K} <\infty.
\end{equation}
 Moreover, $A\le \|T_a\| \le 4A$. \et

\

\section{Bounds for the approximation numbers}

We recall that, given any $m\in \N$, the $m$--th approximation
number of a bounded operator $T:L^p(\G) \to L^p(\G)$, $a_m(T)$,
is defined by
\[ a_m(T):= \inf \|T-F | L^p(\G) \to L^p(\G)\|,
\]
where the infimum is taken over all bounded linear maps $F:L^p(\G)
\to L^p(\G)$ with rank less than $m$.

A measure of non-compactness of $T$ is given by
\[
\beta (T):= \inf \|T-P | L^p(\G) \to L^p(\G) \|,
\]
where the infimum is taken over all compact linear maps $P:L^p(\G) \to L^p(\G)$.
Since $L^p(\G)$ has the approximation property for $1\le p \le
\infty$, $T$ is compact if and only if $a_m(T) \to 0$ as $m \to
\infty$, and $\beta(T)=\lim_{n\to \infty} a_n(T)$.

\bd
Let $K$ be a subtree of $\G$ and $a\in \G$.
We define:
% \[ J_a(K)\equiv J_a(K,u,v):= \|T_{a,K} |L^p(K)\to L^p(K)\| \]
\[ A(K)\equiv A(K,u,v):=
\left\{
\begin{array}{ll}
 \sup_{f\in L^{p}(K),f\not= 0} \inf_{\alpha \in
\cC} {\|T_{a,K}f-\alpha v~\|_{p,K} \over \|f\|_{p,K}}  & \mbox{ if
} \ \mu(K)>0, \\
0 & \mbox{ if } \ \mu(K)=0,
\end{array}
\right.
\]
where
\[
T_{a,K}f(x):= v(x) \chi_K(x) \int_a^x u(t) f(t) \chi_K(t) dt,
\]
and
\[
\mu(K):= \left\{
\begin{array}{ll}
\int_K|v(t)|^p dt
  & 1\le p <\infty, \\
\int_K|v(t)| dt & p=\infty.
\end{array}         \right.
\]
\ed

\bl
The number $A(K,u,v)$ in Definition 3.1 is independent of $a\in
\G$.
\el
{\bf Proof.}
Denote by $S$ the canonical map of $L^p(K)$ into its quotient by
the space of scalar multiples of $v$. Then $A(K)=\|ST_{a,K}|L^p(K)
\to L^p(K) \|$. For $b\in \G$ we have $T_{b,K}f= v\chi_K \int_b^a fu\chi_Kdt +
T_{a,K}Uf $, where $Uf(t)=-f(t)$ if $t$ lies on the path $(a,b)$
and $f(t)$ otherwise. Clearly $U$ is a linear isometry of $L^p(K)$
onto itself and $ST_{b,K}=ST_{a,K}U$. $\Box $

\bc
For all subtrees $K \subseteq \G $
$$
A(K) \le \inf_{a\in \G} \|T_{a,K}|L^p(K) \to L^p(K)\|.
$$
\ec

\
 Note that if $\Lambda$
is a subtree of $\G$, $a\in \G$ and $b$ the nearest point of
$\Lambda$ to $a$ then $T_{b, \Lambda}= T_{a, \Lambda}$ and
$\|T_{b,\Lambda}\|:= \|T_{b,\Lambda}|L^p(\Lambda)\rightarrow L^p(\Lambda)\|
\le \|T_a \| =: \|T_a|L^p(\G)\rightarrow L^p(\G)\|$. Moreover if
$\Lambda'\subset
\Lambda$, with $c$ the nearest point of $\Lambda'$ to $a$ and
$(b,c)$ a subinterval of an edge of $\Lambda$ then
$
T_{b,\Lambda}f=T_{b,\Lambda}(f\chi_{(b,c)})+T_{c,\Lambda'}(f\chi_{\Lambda'}),$
whence $0\le \|T_{b,\Lambda}\|-\|T_{c,\Lambda'}\|\le
\|u\|_{p',(b,c)}\|v\|_{p,(b,c)}$. This remark yields

\bl
For $1\le p\le\infty, \|T_{x,K}| L^p(K) \to L^p(K) \|$ is continuous in $x$.
\el

Let $x\in \G$. Denote by $\G_{x,i}$ $i=1,\dots ,n_x$ the
non-overlapping subtrees of $\G$ which are the closures of the
connected components of $\G  \setminus \{x\}$, and set $ T_{x,i}
\equiv T_{x,\G_{x,i}}$ and $\|T_{x,i}\| \equiv
\|T_{x,i}|L^p(\G_{x,i}) \rightarrow L^p(\G_{x,i})\|  $ . We suppose that the numbering is done in
the order of descending norms of the $T_{x,i}: L^p(\G_{x,i})\to L^p(\G_{x,i})$.
Clearly $\|T_x\|=\max_{i=1,\dots ,n_x} \|T_{x,i}\|$.

\
Call a point $x\in \G$ {\em {simple}} if there is just one $T_{x,i}$ with
maximal norm, so that $\|T_{x,1}\|>\|T_{x,2}\|.$ If $a$ is a simple
point and $(a,y)$ the first edge of $\G_{a,1}$ then by continuity
either there is a point $z$ of $(a,y)$ which is not simple or $ a
\notin \G_{y,1} $. If the latter, continue the path beginning with
$(a,y)$ along the initial edge of $\G_{y,1}$. By induction thus
define a path $l$ in $\G$ satisfying one of the following:
\begin{description}
\item[(i)] $l$ is finite and its end $b$ is not simple;
\item[(ii)] $l$ is finite, its end $b$ is simple and $\{x: x\succeq_a b\}= \emptyset$;
\item[(iii)] $l$ is infinite.
\end{description}

Now (ii) is impossible since $\lim_{x\to b}
\|T_{x,1}\|=0$, and $\|T_{x,1}\| \ge A(\G)$. Also (iii) implies
$T$ is not compact. For if $x$ is in $l$, $\|T_{x,1}\| \ge A(\G)$
and hence there is a compact subset $K$ of $\G_{x,1}$ and a
function $f$ supported in $K$ with $\|f\|\le 1$ and $ \|T_af\|_K
\ge \frac{1}{2}A(\G)$. It follows that there is a sequence of
disjoint compact sets $K_n$ and functions $f_n$ with the same
property. Then, if $m>n, \|T_a(f_n-f_m)\|_{\G}\ge
\|T_af_n\|_{K_n}\ge\frac{1}{2}A(\G)$. Thus, if $T_a$ is compact,
(i) holds. Moreover, $\|T_b\|=\min_{x\in \G}\|T_x\|.$ For
if $x\not= b$ then $x\notin $
one of $\G_{b,1}$, $\G_{b,2}$, say $\G_{b,2}.$ Then if $\G_{x,j}$
is the subtree containing $b$, $\|T_{x,j} \| \ge \|T_{b,2}\| = \|
T_{b,1}\|. $ From this we have the following result which will be
an important tool for determining a lower bound for $A(K)$ once
Theorem 3.8 below is available.
\bl
Suppose $T_a$ is compact and that there exist $i\neq j$ such that $\|T_{a,i}\|,
\|T_{a,j}\|\le \|T_a\|$. Then
\[ \min\{\|T_{a,i}\|,
\|T_{a,j}\|\} \le \min_{x\in \G} \|T_x\| .
\]
\el
{\bf Proof.} The result is a consequence of the discussion preceding the
lemma if $a$ is not simple. If $a$ is simple, then,
with $b$ the non-simple end-point of the path $l$ in (i) above,
and $\min \{\|T_{a,i}\|, \|T_{a,j}\| \} = \|T_{a,i}\| $, say, we
have $\|T_{a,i}\| < \|T_a\| $ and $\G_{a,i}$ is a subtree of some tree
$\G_{b,k}$. Thus
\[
\|T_{a,i}\| \le \|T_{b,k}\| \le \|T_b\| = \min_{x\in \G}\|T_x\|.
\]
{$\Box$}

In the next two lemmas $\|\cdot \|_{p, \mu} $ denotes the norm in
$L^p(\G, d\mu)$, where $d\mu (t) = |v(t)|^pdt $.

\bl
 If $1<p\le \infty$ there is a unique scalar $c_f$ such
that $\|f-c_fe\|_{p,\mu}= \inf_{c\in \cC } \|f-ce\|_{p,\mu}$ for $e \neq 0,
e\in L^p(\G,d\mu)$.
\el
{\bf Proof.} Since $\|f-ce\|_{p,\mu} $ is continuous in $c$ and tends to $\infty$ as
$c \rightarrow \infty $, the existence of $c_f$ is guaranteed by the local
compactness of $ \cC $. For $1<p<\infty$ the uniqueness follows from
the uniform convexity of $L^p(\G,d\mu)$. Let $p=\infty$, and suppose that there are two
values of $c_f, c_1\neq c_2$. This yields the contradiction
$\|f-(1/2)(c_1+c_2)\|_{p,\mu}<\|f-c_1\|_{p,\mu}$.{$\Box$}

\bl The map $f \to c_f : L^p(\G,d\mu) \rightarrow \cC $ is continuous for $1<p\le \infty$.
\el
{\bf Proof.} Suppose that $c_{g_n} \rightarrow c $ as $g_n \rightarrow f $.
Then
\[
\|g_n-c_f\|_{p,\mu} \ge \|g_n-c_{g_n}\|_{p,\mu}
\]
and so
\[
\|f-c_f\|_{p,\mu} \ge \|f-c\|_{p,\mu}
\]
which gives $c=c_f $
{$\Box$}

\bt Let $1<p\le \infty$. If $T_a$ is compact
$A(\G)=\min_{x\in \G}\|T_x|L^p(\G)\to L^p(\G)\|$.
\et
{\bf Proof.} There is a non-simple point $b$ at
which $\|T_x\|$ attains its minimum. If $\alpha < \|T_b\|$ there
exist $f_i, i=1,2$, supported in $\G_{b,i}$ with $\|f_i\|=1,
\|T_bf_i\|>\alpha $, and $f_1$ positive, $f_2$ negative. Clearly
the same is true of the corresponding values of $c_f$, say $c_1,
c_2$. Then by continuity there is a $\lambda \in [0,1]$ such that
$c_g=0$ for $g=\lambda f_1 + (1-\lambda) f_2$, and $\|T_b
g\|^p=\lambda^p\|T_bf_1\|^p + (1-\lambda)^p \|T_bf_2\|^p
>\alpha^p \|g\|^p$. Then, by Lemma 3.6,
\[
 A(\G) \ge \inf_c \|(T_b-cv)g\|/\|g\| = \|T_bg\|/\|g\| > \alpha.
\]
Since $\alpha <\|T_b\| $ is arbitrary, $A(\G) \ge \|T_b\|$  and the result
follows from Corollary 3.3.{$\Box$}

The next lemma establishes an important geometrical property of a
tree which is an essential ingredient of the subsequent analysis.
First we make some observations.

Suppose $w$ is a non--negative function defined on the set of
all closed subtrees of a tree $\G$, satisfying
\begin{equation} \label{w}
X\subseteq  Y \Rightarrow  w(X) \le w(Y).
\end{equation}
Define
\begin{itemize}
\item

\[
N_{\varepsilon}(\G)=\min_{\cF \in \cS_{\varepsilon}(\G)}
\# \cF
\]
where  $\cS_{\varepsilon}:= \{\cF; \cF \mbox{ is a set of
non--overlapping closed subtrees of $\G$ such that }$ $i)
\cup_{X\in \cF}X=\G, ii) X\in \cF \Rightarrow  w(X) \le
\varepsilon \}$;
\item
\[
M_{\varepsilon}(\G)=\max_{\Ge \in \cL_{\varepsilon}(\G)}
\# \Ge
\]
where $\cL_{\varepsilon}:=\{\Ge; \Ge \mbox{ is a set of
non--overlapping closed subtrees of $\G$ such that }$
 i) $\cup_{X\in \Ge}X=\G$ \, ii) $\#\{X; X\in \Ge, w(X)\le
\varepsilon\} \le 1 \} $
\end{itemize}

Two non--overlapping closed subtrees of $\G$ can have at most
one point in common, for otherwise $\G$ would contain a cycle. A
{\em chain} ${\cal C}$ of closed subtrees is a sequence $X_1, \dots ,X_l$ of
closed subtrees such that $X_i \cap X_{i+1}=\{x_i\}$ $(i=1, \dots
,l-1)$ where the $x_i$ are distinct. The {\em length} of ${\cal C}$ is $l$.

There is a set $\cF \in
\cS_{\varepsilon}(\G)$ with $\#\cF=N_{\varepsilon}(\G)$
(possibly $\infty$). Let ${\cal C}$ be a chain of elements of $\cF$ of
maximal length $l$. Then we have the following:

\begin{description}
\item[(i)]
If $l=1$ then $\#\cF=1$ and so $w(\G)\le \varepsilon$.

\item[(ii)]
If $l=2$ define $Y=X_1\cup X_2$. Then if $\G \not= Y$,
$\G \setminus Y^o$ is a closed subtree and $N_{\varepsilon}(\G
\setminus Y^o)= N_{\varepsilon}(\G)-2$. Moreover $Y$ is a closed
subtree of $\G$ and $w(Y) > \varepsilon$. For since $l$ is maximal
$x_1$ lies in every member of $\cF$ so that $\G \setminus Y^o$ is
a closed subtree. Also $\G \setminus Y^o=\cup_{X\in \cF'} X$,
where $\cF'=\cF \setminus \{X_1,X_2\}$, which implies that
$N_{\varepsilon}(\G \setminus Y^o) \le N_{\varepsilon}(\G)-2.$ If
$N_{\varepsilon}(\G \setminus Y^o) < N_{\varepsilon}(\G)- 2$, there
exists $\cF''$, a suitable covering of $\G \setminus Y^o$ with $\#\cF''=
N_{\varepsilon}(\G \setminus Y^o)$, and then $\cF''\cup
\{X_1,X_2\}\in \cS_{\varepsilon}(\G)$ which contradicts the
definition of $N_{\varepsilon}(\G)$. Finally $w(Y)>\varepsilon $
for if not, on taking $\cF'''= \cF' \cup \{Y\}$ we have a contradiction.

\item[(iii)]
If $l \ge 3,$ suppose $Z_1, \dots ,Z_k$ are the sets
in $\cF$ which meet $X_{l-1}$ in a point different from
$x_{l-2}$. Then since $l$ is maximal $X_1, X_2, \dots ,X_{l-1},
Z_i$, is a chain of maximal length for $i=1, 2, \dots ,k$. Either
a) $k=1$ or b) $k>1$. If a), we have $Z_i=X_l$ and we take
$Y=X_l\cup X_{l-1}$, so $Y$ is a closed subtree. Then $\G \setminus Y^o $ is a closed
subtree of $\G$. Moreover $N_{\varepsilon}(\G \setminus
Y^o)=N_{\varepsilon}(\G)-2$, $w(Y)>\varepsilon$ by an argument
similar to that of (i). If b), define $a_i$ by $\{a_i\}=X_{\l-1}\cap Z_i$
and, for $i\neq
j, a_{ij}=a_i \wedge a_j=\max \{u: u \preceq a_i \mbox{ and }
u\preceq a_j\}$ in the ordering of $\G$ arising from taking
$x_1$ as root. Then $a_{ij}\in X_{l-1}.$ Define $\varrho_i:=\#\{a_{jk};
a_{jk} \preceq a_i \}. $ Without loss of generality we may
suppose $\varrho_1=\max_i \varrho_i $ and that $a_{jk} \preceq
a_{12}$ for all $a_{jk} \preceq a_1$. Define $Y=\{x;
(\exists u) a_{12} \prec u \preceq a_1 \mbox{ or }
a_{12}\prec u \preceq a_2 \mbox{ and } u \preceq x \}
\cup \{a_{12} \} $. Then $Y$ is a closed subtree of $\G$ and
$Y\subseteq Z_1 \cup Z_2 \cup X_{l-1}.$ If $\G'=\G \setminus Y^o$,
$\G'$ is a closed subtree of $\G$ and it follows by arguments
similar to those in (ii) that, since ${\cal C}$ is a chain of maximal length in
$\cF$ and $\cF$ is a covering of $\G$, $w(Y)>\varepsilon $ and
\begin{equation} \label{*} N_{\varepsilon}(\G)-3 \le N_{\varepsilon}(\G') \le
N_{\varepsilon}(\G) -2. \end{equation}
\end{description}

\bl If $N_{\varepsilon}(\G) <\infty$ then $M_{\varepsilon}(\G) \ge
{1\over 3} N_{\varepsilon}(\G)$.
\el
{\bf Proof.} The proof is by
induction on $N_{\varepsilon}(\G)$. If $N_{\varepsilon}(\G)=1,
w(\G) \le \varepsilon$ and $M_{\varepsilon}(\G)=1$. If
$N_{\varepsilon}(\G)=2, l=2$ and $\G=X_1 \cup X_2$ and
$M_{\varepsilon}(\G) \ge 1$. If $N_{\varepsilon}(\G)>2$ the
arguments of i) ii) iii) above and induction prove the result. For
by i) $l \ne 1$ and by ii) or iii) $N_{\varepsilon}(\G \setminus
Y^o)<N_{\varepsilon}(\G)$ so we may suppose $M_{\varepsilon}(\G
\setminus Y^o)\ge {1 /3} N_{\varepsilon}(\G \setminus Y^o)$ and
then there exists $\Ge' \in \cL_{\varepsilon}(\G \setminus Y^o)$
with $\#\Ge'=M_{\varepsilon}(\G \setminus Y^o)$. But then if
$\Ge=\Ge'\cup \{Y\}$, $\Ge \in {\cL}_{\varepsilon}(\G)$ since $\omega(Y)>
\varepsilon $ and
$M_{\varepsilon}(\G) \ge \#\Ge = M_{\varepsilon}(\G \setminus
Y^o)+1\ge {1 \over 3} N_{\varepsilon}(\G \setminus Y^o)+1 \ge
{1\over 3}N_{\varepsilon}(\G).$ {$\Box$}

\bl Let $K$ be a compact tree, let $w$ be a function satisfying
(\ref{w}) and, for every $(c,d) \in E(K)$, suppose that $w(c,.)$ is
a continuous function on $(c,d)$. Then there exists a set ${\cal G}$ of
non-overlapping subtrees of $K$ with $\omega(X)\ge \varepsilon $
for $ X\in {\cal G}$ and
\begin{equation} \label{**}
 \# {\cal G} \ge
 N_{\varepsilon}(K) - 3\# E(K).
\end{equation}
\el
{\bf Proof.} Let $\cF \in S_{\varepsilon}(K)$ and $\#\cF =
N_{\varepsilon}(K)$. For $(c,d)\in E(K)$ define ${\cal K} := \{X: X\in
\cF,|X\cap(c,d)|>0\}$ to be such that $\# {\cal K} >3$.  Set ${\cal
K}
=\{X_1,\dots,X_n\}$, where $n=\#{\cal K}$ and $X_1 \preceq X_2 \dots
\preceq X_n$, where $X\preceq Y $ means that $x\preceq y$ for all
$x\in X,y\in Y$. Then $X_i=(a_i,b_i)\subset (c,d)$ for $1<i<n$. By
the continuity of $\omega $ on $(c,d)$ and the minimality of
$\#\cF$, there exists $b_2^{\prime} \in(b_2,b_3) $ such that
$\omega(a_2,b_2^{\prime})=\varepsilon $. It follows that there are
non-overlapping sub-intervals $X_2^{\prime},\dots,X_k^{\prime} $ of
$(c,d)$, where $k=n-2$, for which
$\omega(X_j^{\prime})=\varepsilon$. The lemma follows from this.{$\Box$}

\bd
Let $K$ be a subtree $\G$. We denote by $\cP(K)$ the set of
partitions $\{\G_i: i=1, \dots ,n\}$ of $K$ (i.e $\cup_{i=1}^{n}
\G_i=K$) by subtrees $\G_i$ of $K$
 such that $|\G_i\cap \G_j|=0$ for $i\not= j$. For a
given $\varepsilon >0$ we define
\begin{eqnarray*}
N(K,\varepsilon)\equiv  N(K,\varepsilon,u,v):= \min  \{n:
\exists & \{\G_i: i=1, \dots ,n\}\in \cP(K)&\\
& \mbox{ and }
A(\G_i,u,v)\le \varepsilon\};&
\end{eqnarray*}
\begin{eqnarray*}
M(K,\varepsilon) \equiv  M(K,\varepsilon,u,v):= & \max \{m :
\exists \mbox{ non-overlapping subtrees } \G_i \subseteq K,&\\
&i=1,\dots,m, \mbox{ such that }
A(\G_i,u,v)> \varepsilon \}.&
\end{eqnarray*}
\ed
Note that, with $\omega (\cdot) = A(\cdot)$,
$M_{\varepsilon}(\G)\le M(\G, \varepsilon)+1$.

Hereafter we shall write $T, T_K $ for $T_a, T_{a,K}$
respectively, unless there is a possibility of confusion. The
following lemma will yield a one-dimensional approximation to $T$
on a subtree $K$ of $\G$. We recall the notation
\[
\mu(K) = \int_K d\mu,\;\;\; d\mu = v^p dx.
\]

\bl
 Let $K$ be a subtree of $\G$ and $v \in L^p(K)$, $1\le p \le
\infty$, with $\mu(K) \not= 0$. Then there exists $w_K\in \{L^p(K;
d\mu)\}^*$ (in case $p=\infty$, $w_K$ is also from
$\{L^{\infty}(K)\}^*$ ) such that:
\[w_K(1)=1,\]
\[\|w_K\|_{\{L^p(K,d\mu)\}^*} = {1\over \|v\|_{p,K}},
(\| w_K\|_{\{L^{\infty}(K)\}^*}=1 \mbox{ for }p=\infty )
\]
and
\begin{equation}\inf_{c\in \cC}\|(\varphi-c)v\|_{p,K} \le \|(\varphi -
w_K(\varphi))v\|_{p,K} \le 2 \inf_{c\in \cC}
\|(\varphi-c)v\|_{p,K}
\end{equation}
for all $\varphi \in L^p(K,d\mu)$. In the case $p=2$
\begin{equation}
\inf_{c\in \cC}\|(\varphi-c)v\|_{L^p(K)}= \|(\varphi -
w_K(\varphi))v\|_{L^p(K)}
\end{equation}
where
\[
w_K(\varphi) = \frac{1}{\mu(K)}\int_K \varphi d\mu.
\]
\el
{\bf Proof.} Define the linear function $w$ on the constants
in $L^p(K,d\mu)$ with $w(c.1)=c$. Then $w(1)=1$ and
$\|w\|^p=1/{\mu(K)}$ for $1\le p <\infty$ while $\|w\|=1$ when
$p=\infty$. The existence of $w_K$ follows by the Hahn-Banach theorem, and
(3.4) is immediate. The case $p=2$ is obvious.
{$\Box$}

\br
For $1\le p < \infty$ and $\G=(0,\infty)$ the choice
$w_K(\varphi)={1\over \mu(K)} \int_K \varphi d\mu$ is
appropriate, and was that used in [2]. In [6] Lemma 2.4,
when $p=\infty$ and $\G=(a,b)$, $w_K$ was defined as the limit
along a filter base of subsets $w_{\beta}$ of the unit ball in
$\{L^{\infty}(K)\}^*$ defined by
\[w_{\gamma}(\varphi):={1\over |A_{\gamma}|} \int_{A_{\gamma}}
\varphi(x) dx, \qquad \varphi\in L^{\infty}(K). \]
\er

\bl
Let $K$ be a subtree of $\G, a\in K, 1<p\le \infty$, and suppose that
$T_{a,K}$ is compact. Then
\[
A(K)=\|T_{a,K}-v \varpi |L^p(K)\rightarrow L^p(K)\|,
\]
where $ \varpi $ is the bounded linear functional
\[
\varpi(f)= \int_a^bfu
\]
and $b\in K$ is such that $A(K)=\|T_{b,K}|L^p(K)\rightarrow
L^p(K)\|$.
\el
{\bf Proof.}
We know from Theorem 3.8 that there exists a $b\in K$ such that
\[
A(K,u,v)=\|T_{b,K}|L^p(K)\rightarrow L^p(K)\|= \|T_{b,K}U|L^p(K)\rightarrow
L^p(K)\|,
\]
where $U$ is the linear isometry defined in the proof of Lemma
3.2, and with respect to which
\[
T_{b,K}f=v\chi_K\int_b^afu\chi_K +T_{a,K}Uf.
\]
Thus, if we define
\[
Pf(x)=v(x)\chi_K(x)\int_a^bUfu\chi_K,
\]
we have
\[
\|T_{a,K}-P|L^p(K)\rightarrow L^p(K)\| = A(K,u,v),
\]
and the lemma follows.{$\Box$}

\bl Let $\varepsilon >0$, $1\le p \le \infty$ and let $K$ be a
subtree of $\G$. If $N:=N(K,\varepsilon,u,v)<\infty$ then
\[
 a_{N+1}(T_K)\le \gamma_p \varepsilon,
\]
where $\gamma_p =2$ when $p\neq 2$ and $\gamma_2=1$.
\el
{\bf Proof.} Let $\{\G_i\}_1^N$ be the partition of $K$ which defines $
N := N(K,\varepsilon,v,u)$ in Definition 3.11 and set
$Pf=\sum_{i=1}^N P_if$ where
\[
P_if(x):= \chi_{\G_i}(x) v(x) \left[ \int_a^{a_i} uf
 + w_i\left( \int_{a_i}^x u f \chi_{\G_i}
\right) \right],
\]
$T_K\equiv T_{a,K}, a_i$ is the point in $\G_i$ nearest $a$ and
$w_i \equiv w_{\G_i}$ is the linear function from Lemma 3.12.
Then $\rank P \le N$ and, on using Lemma 3.12, and setting
$T_i\equiv T_{\G_i} $,

\begin{eqnarray*}
\|(T_K-P)f\|^p_{p,K}&= & \sum_{i=1}^N \|(T_K-Pf) \|^p_{p,\G_i}
\\
 &= & \sum_{i=1}^N
\|T_i-w_i\left(\int_{a_i}^{\cdot} uf\chi_{\G_i} dt \right) v(\cdot)
\|^p_{p,\G_i}  \\
 &\le & \gamma_p^p \left( \max_{i=1,\dots ,N}
A(\G_i,u,v) \right)^p \|f\|^p_{p,K}\\
 &\le & (\gamma_p\varepsilon)^p \|f\|^p_{p,K},
\end{eqnarray*}
whence the lemma.{$\Box$}

\bl
Let $ 1<p\le \infty, \varepsilon >0,$ and suppose that $T_K$ is
compact. Then, with $N=N(K,\varepsilon,u,v)<\infty$,
\[
a_{N+1}(T_K) \le \varepsilon.
\]
\el
{\bf Proof.}
The same argument as in Lemma 3.15 applies but with
\[
P_if(x)=\chi_{\G_i}(x)v(x)\left[\int_a^{a_i}uf + \varpi_i
\left(\int_{a_i}^{b_i}uf\chi_{\G_i}\right)\right],
\]
where $ \varpi_i, b_i$ are as in Lemma 3.14 corresponding to
$K=\G_i$.{$\Box$}

\bl Let $\varepsilon >0$, $1\le p \le \infty$ and let $K$ be a
subtree of $\G$. Let $\{ \G_i; i=1,
\dots ,M\}$ be a set of non-overlapping subtrees of $K$ such that
$A(\G_i)\ge \varepsilon $ for all $1\le i \le M$. Then
\[
a_M(T_K) \ge \varepsilon.
\]
\el {\bf Proof.} Let $\lambda \in (0,1).$ From the definition of
$A(\G_i)$ we know that for $i=1, \dots ,M$ there exist $\phi_i \in
L^p(\G), \|\phi_i \|_p=1$, with support in $\G_i $ such that
\[
\inf_{\alpha \in \cC} \|T_K\phi_i-\alpha v\|_{p,\G_i}
>\lambda A(\G_i) \ge \lambda \varepsilon.
\]

Let $P:L^p(K) \to L^p(K)$ be bounded with $\rank (P)<M$. Then,
there are constants $\lambda_1, \dots , \lambda_M$ not all zero,
such that
\[ P\phi =0, \qquad \qquad \phi:=\sum_{i=1}^M
\lambda_i \phi_i.
\]
Then, noting that the following summation is over $\lambda_i
\neq 0$, and denoting by $a_i$ the point of $\G_i$ nearest $a$,
where $T_K =T_{a,K}$,
\begin{eqnarray*}
a_M^p\|\phi\|_{p,K}^p & \ge & \|T_K\phi -  P\phi \|^p_{p,K} =  \|T_K\phi
\|^p_{p,K} = \sum_{i=1}^M \| (T_K\phi)\chi_{\G_i} \|^p_{p,K}
\\
 & = & \sum_{i=1}^M \| \chi_{\G_i}(x)v(x) \left( \int_{a_i}^x
\lambda_i \phi_i(t) u(t) \chi_{\G_i}(t) dt \right.  \\ &&\qquad
\qquad \left.  + \int_a^{a_i} \phi(t) u(t) dt \right)
\|^p_{p,K} \\
 &= & \sum_{i=1}^M
 \|\left(T_K\phi_i(x)+v(x) {\eta_i \over \lambda_i} \right) \lambda_i
 \|^p_{p,\G_i},\\
 & & \mbox{ where }   \eta_i:=
\int_a^{a_i} \phi(t) u(t) dt, \\
&\ge & \sum_{i=1}^M
  \inf_{\alpha \in \cC} \|( T_K \phi_i(x) - v(x) \alpha)
 \|^p_{p,\G_i} |\lambda_i |^p \\ &\ge &(\lambda
\varepsilon)^p \sum_{i=1}^M \| \phi_i \|^p_{p,\G_i}
|\lambda_i |^p  \ge (\lambda \varepsilon)^p \|\phi \|^p_{p,K}.
\end{eqnarray*}
{$\Box$}

\bt
Let $1\le p \le \infty, \varepsilon>0 $ and
$N=N(\G,\varepsilon,u,v)<\infty $(see Definition 3.11).
Then
\[
a_{N+1}(T) \le \gamma_p \varepsilon
\]
where $\gamma_p = 2 $ when $p\neq 2$ and $\gamma_2 =1$, and
\[
a_{[{N \over 3}] -1}(T) >\varepsilon.
\]

The measure of non-compactness of $T, \beta(T)$, satisfies
\[
\beta(T) :=\lim_{n \rightarrow \infty}a_n(T) \asymp
\inf\{\varepsilon: N(\G,\varepsilon,u,v)<\infty\},
\]
where the symbol $\asymp$ means that the quotient of the two sides
lies between positive constants. Hence, $T$ is compact if and only if
$N(\G,\varepsilon,u,v)<\infty$ for all $ \varepsilon >0$. If $T$ is
compact, $\gamma_p=1$ for $1<p\le\infty$.
\et
{\bf Proof.} The first
inequality follows from  Lemma 3.15. The second inequality follows from Lemmas
3.9 and 3.17 on putting $w(\cdot):= A(\cdot)$. The two inequalities
imply the result about the measure of non-compactness, and hence
the compactness, of $T$. The last statement is a consequence of
Lemma 3.16.{$\Box$}

From Lemmas 3.10, 3.16, 3.17 and 3.18 with $w(X):=A(X)$ we derive
\bl Let
$1< p < \infty$, $\varepsilon >0$,
 and let $K$ be a compact tree. Then
\[a_{N+1}(T_K) \le  \varepsilon\]
and
\[a_{M-1}(T_K) >  \varepsilon,
\]
where $N=N(K,\varepsilon,u,v),and $ $M\ge N-3 \# E(K).$
\el

Since, by Lemma 3.2, $A(\G_a)$ is independent of $a\in \G$, the
approximation numbers are independent of $a\in \G$. Note that the
above proof of Lemma 3.19 requires $A(c,.)$ to be continuous (see
Lemma 3.10). This is not true for $p=1$ or $\infty$.

\section{Local properties of $A$}

In this section we establish properties of the function $A$
in Definition 3.1 which will be needed for the local
asymptotic results in the next section.

\bl Let $u, v$ be constants over a real interval $I=(a_1,b_1)$ and
$1\le p \le \infty$. Then $A(I,u,v)=|v| |u| |I| \alpha_p$, where
$\alpha_p=A((0,1),1,1)$.
\el {\bf Proof.} We have

\begin{eqnarray*}
A(I,u,v) &= & \sup_{\|f\|_{p,I}\le 1} \inf_{\alpha \in \cC} \|v
\left(\int_{a_1}^x uf(t) dt -\alpha \right) \|_{p,I} \\
&= & |v| |u| \sup_{\|f\|_{p,I}\le 1} \inf_{\alpha \in \cC} \|
\int_{a_1}^x f(t) dt -\alpha \|_{p,I} \\
&= & |v| |u| |I| \sup_{\|f\|_{p,(0,1)} \le 1} \inf_{\alpha \in
\cC} \|\int_0^x f(t) dt -\alpha \|_{p,(0,1)} \\
&= & |v| |u| |I| A((0,1), 1, 1).
\end{eqnarray*}
{$\Box$}

Note that $ \alpha_1= \alpha_{\infty}=1/2$ and $\alpha_2=1/\pi $.
\bl
Let $1\le p\le \infty$, $u_1, u_2 \in
L^{p'}(\G)$ and $v\in L^p(\G)$. Then
\[|A(\G,u_1,v)-A(\G,u_2,v)| \le \|v\|_{p,\G} \|u_1-u_2\|_{p',\G}.
\]
\el {\bf Proof.} We have
\begin{eqnarray*}
 A(\G,u_1,v) &=& \sup_{\|f\|_p =1} \inf_{\alpha \in \cC}
\|v(x)\left( \int_a^x (u_1(t)-u_2(t)+u_2(t)) f(t) dt -\alpha
\right) \|_p \\
 & \le & \sup_{\|f\|_p=1} \inf_{\alpha \in \cC}
\left[ \| v(x) \int_a^x (u_1(t)-u_2(t))f(t)dt \|_{p,\G}  \right.
\\ &+&  \qquad \qquad \qquad \left. \|v(x)\left( \int_a^x u_2(t)
f(t)dt -\alpha \right)\|_p \right]
\\ & \le &
\sup_{\|f\|_p=1} \inf_{\alpha \in \cC} \left[ \|v\|_p
\|u_1-u_2\|_{p'}+ \|v(x)\left( \int_{a_1}^x u_2(t) f(t) dt
-\alpha \right) \|_p \right] \\
&\le &\|v\|_p
\|u_1-u_2\|_{p'} +A(\G,u_2,v).
\end{eqnarray*}
The same holds with $u_1,u_2$ interchanged.
 {$\Box$}

\bl
Let $1\le p \le \infty$, $u\in L^{p'}(\G)$, and
$v_1, v_2 \in L^p(\G)$. Then
\[|A(\G,u,v_1) -A(\G,u,v_2)| \le 2 \|v_1-v_2\|_p
\|u\|_{p'}.
\]
\el

{\bf Proof.} We have
\begin{eqnarray*}
A(\G,u,v_1) &=& \sup_{\|f\|_p=1} \inf_{\alpha \in \cC}
\|v_1(x) \left[ \int_a^x u(t) f(t) dt-\alpha \right] \|_p \\
&= & \sup_{\|f\|_p=1} \inf_{|\alpha |\le \|u\|_{p'}
\|f\|_p} \|v_1(x) \left[ \int_a^x u(t) f(t) dt - \alpha
\right] \|_p.
\end{eqnarray*}
Since
\[
\|v_1(x)[\int_a^xu(t)f(t)dt -\alpha]\|_p
\le \|v_1-v_2\|_p\|u\|_{p'}\|f\|_p + \|(v_1-v_2)\alpha\|_p
+\|v_2[\int_a^x u(t)f(t)dt -\alpha]\|_p
\]
it follows that
\begin{eqnarray*}
A(\G,u,v_1) &\le &   2\|v_1 - v_2\|_p
\|u\|_{p'} + \sup_{\|f\|_p=1} \inf_{|\alpha |\le
\|u\|_{p'}} \|v_2(x) \left[ \int_a^x u(t) f(t) dt -\alpha
\right] \|_p \\
 &= & 2\|v_1-v_2 \|_p \|u\|_{p'}+A(\G,u,v_2).
\end{eqnarray*}
Similarly with $v_1$ and $v_2$ interchanged. {$\Box$}

\bl Let $1< p < \infty$ and suppose that  $K_1,K_2,
K_1 \subset K_2$, are compact subtrees of $\G$. Then
\[ |A(K_1) - A(K_2)| \to 0 \qquad \mbox{ as } |K_2 \setminus
K_1| \to 0 \]
\el

{\bf Proof.} We see that $A(K_1)=A(\G,u\chi_{K_1},v\chi_{K_1})$
and $A(K_2)=A(\G,u\chi_{K_2},v\chi_{K_2}).$  The lemma then
follows from Lemmas 4.2 and 4.3.{$\Box$}

In order to treat the cases $p=1,\infty$ we need the following
terminology.

\bd
Let $u\in L^{\infty}(\G).$ Then
\[u_s:= \lim_{\varepsilon \to 0_+} \|u
\chi_{B(t,\varepsilon)}\|_{\infty,\G}
\]
where $B(t,\varepsilon
)$ is the ball center $t$, radius $\varepsilon $ on $\G$.
\ed

\bd Let $g$ be a function defined on a real interval $I$. Then
\[g^*(x):= \inf\{t;g_*(t) \ge x\}, \]
where $g_*(t):=|\{x\in I; g(x) \ge t \}|.$
The function $g^*$ is the non--increasing rearrangement of
$g$.
\ed

Note that since we have $\ge$ in the definitions above, $g_*$ and
$g^*$ are left--continuous functions. For the case $p=\infty$ we
have the following two lemmas.

\bl Let $I$ be a bounded interval, $\gamma, \delta \in \R $ with
$\delta \ge v_s(t) \ge 0$ on $I$, and let $p=\infty$. Then
\[A(I;\gamma, \delta) \ge A(I; \gamma,v_s) \ge 1/2 |\gamma|
\|(v_s \chi_I)^*(t) t \|_{\infty,(0,|I|)}. \]
\el

{\bf Proof.} See [6; Lemma 4.5]. {$\Box$}

\bl Let $I$ be a bounded interval, $\gamma, \delta \in \R $ with
$\delta \ge v_s(t) \ge 0$ on $I$ and let $p=\infty$. Then for
any $\alpha >1$
\[ A(I;\gamma, \delta)-A(I;\gamma, v_s) \le {\alpha \over 2}
\int |\gamma|(\delta -v_s(t))dt +{|\gamma| \delta |I| \over 2
\alpha }. \]
\el

{\bf Proof.} See [6, Lemma 4.6]. {$\Box$}

For the case $p=1$ we have

\bl Let $I$ be a bounded interval, $\gamma, \delta \in \R $ with
$\delta \ge u_s(t) \ge 0$ and $\infty > u_1(t) \ge u_2(t) \ge 0$
on $I$. Then for $p=1$
\[A(I;u_1, \gamma) \ge A(I; u_2, \gamma) \]
and
\[A(I;u_s, \gamma) \ge 1/4 |\gamma|
\|(u_s \chi_I)^*(t) t \|_{\infty,(0,|I|)}. \]
\el

{\bf Proof.} In the first inequality
\[ A(I;u_1, \gamma)= |\gamma | \sup_{\|f\|_{1,I}\le 1}
\inf_{\alpha \in \cC} {\| \int_a^x u_1f -\alpha \|_{1,I}  \over
\|f \|_{1,I}}.
\]
For any $\|f_2\|_{1,I} \le 1$ there exists $f_1$
such that $\|f_1\|_{1,I} \le \|f_2\|_{1,I} \le 1$ and
$\int_a^xu_1f_1=\int_a^xu_2f_2$. (Put $f_1(t):=f_2(t){u_2(t) /
u_1(t)}$ if $u_1(t)\not= 0$ and $f_1(t):=f_2(t)$ otherwise.) Then
$A(I;u_1,\gamma) \ge A(I;u_2,\gamma) \ge 0$.

In the second inequality, we have
\begin{eqnarray*}
A(I;u_s,\gamma)& = & |\gamma| \sup_{\|f\|_{1,I}=1} \inf_{\alpha
\in \cC} \|\int_a^x u_s f -\alpha \|_{1,I} \\
& \ge & |\gamma| \sup_{\|f\|_{1,I}=1} \inf_{\alpha \in \cC} \|
\int_a^x \chi_{M_{\beta}} \beta f - \alpha \|_{1,I}
\end{eqnarray*}
where $M_{\beta}:=\{y\in I; u_s(y) \ge \beta\}$ and $0 \le \beta
\le \|u_s\|_{\infty,I}$.
Put $f=\delta_{x_{\beta}}$, where $x_{\beta} \in I =(a,b)$ and
$|M_{\beta} \cap (a,x_{\beta})|=|M_{\beta} \cap
(x_{\beta},b)|={1/2} |M_{\beta}|.$
Then
\begin{eqnarray*}
A(I; u_s, \gamma) & \ge & |\gamma|\inf_{\alpha \in \cC}
\|\chi_{(x_{\beta},b)} \beta -\alpha \|_{1,I} \\
& = & |\gamma| \inf_{\beta > \alpha > 0} (\alpha
(x_{\beta}-a), (\beta - \alpha)(b-x_{\beta})) \\
& = & |\gamma| |\beta| \inf_{1 > \alpha >0} (\alpha
(x_{\beta}-a), (1 - \alpha)(b-x_{\beta})) \\
& = & |\gamma| |\beta| {b-x_{\beta} \over b-a} (x_{\beta} -a) \\
& \ge & |\gamma| |\beta| |{M_{\beta} \over 2 }|
{b-a-|M_{\beta}/2| \over b-a} \\
& \ge & |\gamma| |\beta| |M_{\beta}| {1 \over 2}
{b-a-({b-a \over 2}) \over b-a} \\
& = &{1 \over 4} |\gamma| |\beta| |M_{\beta}|.
\end{eqnarray*}

Hence, we have for every $0 \le \beta \le \|u_s\|_{\infty,I}$
\[A(I; u_s,\gamma) \ge |\gamma| |\beta| |M_{\beta}| {1 \over
4} \]
and so
\[A(I; u_s,\gamma) \ge |\gamma| {1 \over 4} \|t (u_s
\chi_I)^*(t) \|_{\infty,(0,|I|)}.
\]
{$\Box$}

\bl
Let $I$ be a bounded interval, $\gamma, \delta \in \R $ with
$\delta \ge u_s(t) \ge 0$ on $I$ and $p=1$. Then for any $\alpha
>1$
\[A(I;\delta, \gamma)-A(I;u_s, \gamma) \le {\alpha \over 2}
\int_I |\gamma| (\delta - u_s(t)) dt + {|\gamma| \delta |I| \over
2 \alpha }.
\]
\el {\bf Proof.} From Lemmas 4.9 and 4.1 we have
\[ 0 \le A(I;\delta, \gamma)-A(I; u_s,\gamma) \le {1 \over 2}
|\gamma| |\delta| |I| -{1 \over 4} |\gamma| \|(u_s \chi_{I})^*(t)
t \|_{\infty}.
\] The rest of the proof is similar to that in [6,
Lemma 4.6] on using Lemma 4.8 instead of [6, Lemma 4.5]. {$\Box$}

\section {Local asymptotic results}

{ {\bf 5.1 Case $1 < p < \infty$} }

\bl Let $K$ be a compact subtree of $\G$ and $1<p<\infty$. Then
\[ \alpha_p \int_K |u| |v| = \lim_{\varepsilon \to 0_+}
\varepsilon N(K, \varepsilon, u, v),\]
\[ \alpha_p \int_K |u| |v| = \lim_{\varepsilon \to 0_+}
\varepsilon M(K, \varepsilon, u, v),
\]
where $N(K, \varepsilon, u, v)$ and $ M(K, \varepsilon, u, v)$ are
defined in Definition 3.11, and $\alpha_p = A((0,1),1,1)$ (see
Lemma 4.1).
\el {\bf Proof.} Since $K$ is a compact tree it has a bounded
number of vertices, i.e. $K$ is a finite union of intervals. The
argument in [2, Theorem 5], with A replacing the function $L$
there, continues to go through and yields the first equation of the lemma.
The second identity follows from the first identity and Lemma 3.10 since
$A$ is a continuous function on an interval. {$\Box$}

\bl
Let $1<p<\infty$. Then
\[
\alpha_p \int_{\G} |u| |v| \le \liminf_{\varepsilon \to 0+}
\varepsilon N(\G, \varepsilon,u,v).
\]
\el
{\bf Proof.} There exist
compact subtrees $\G_n$ of $\G$, $n=1,2,\dots $ such that $\G_n
\subset \G $ and $\G_n \to \G$ (i.e. $|\G-\G_n|\rightarrow 0 $) as
$n\to \infty$. By Lemma 5.1 we have
\[
\alpha_p \int_{\G_n} |u| |v| = \lim_{\varepsilon \to 0_+}
\varepsilon N(\G_n, \varepsilon,u,v) \le \liminf_{\varepsilon \to
0_+} \varepsilon N(\G, \varepsilon,u,v),
\]
whence the result.
{$\Box$}

\bt Let $1<p<\infty$, $u \in L^{p'}(\G)$ and $v\in L^p(\G)$. Then
\[
 \lim_{\varepsilon \to 0_+} \varepsilon N(\G,\varepsilon,u,v)
 =\alpha_p \int_{\G} |u| |v|,
\]
\[
\lim_{\varepsilon \to 0_+} \varepsilon M(\G,\varepsilon,u,v) =\alpha_p \int_{\G} |u| |v|.
\]
\et
{\bf Proof.}
Let $\cL$ be a maximal such that $\# \cL =M(\G,\varepsilon)\equiv
M(\G,\varepsilon,u,v)$ and $\cS$ a minimal cover such that $\# \cS =
N(\G,\varepsilon)\equiv N(\G,\varepsilon,u,v)$.

Given $\eta>0$, choose a compact subtree $ K \subset \Gamma $ such that
$\|u\chi_{\Gamma \setminus K}\|_{p'} \le \eta $ and
$\|v\chi_{\Gamma \setminus K}\|_{p} \le \eta $.

Set
\begin{eqnarray*}
\cL_1 & = & \{\Gamma':\Gamma' \in \cL, \Gamma' \subseteq K\}, \\
\cL_2 & = & \{\Gamma':\Gamma' \in \cL, \G'\subseteq \Gamma \setminus K \}, \\
\cL_3 & = & \cL\setminus \{\cL_1 \cup \cL_2\},
\end{eqnarray*}
and similarly for $\cS_1, \cS_2, \cS_3$ with respect to $ \cS $.

We know that $\Gamma \setminus K$ is the union of disjoint connected
components
$\{ \Gamma_i^* \}$ and $\# \{ \Gamma_i^* \} \le \# \partial K$. Also, if
$X\in \cL_2 \cup \cS_2$ then $X \subseteq \Gamma_i^*$ for some $i$.
Thus $\# \cL_2 \le \sum_i M(\G_i^*,\varepsilon).$

Consider now the union of $\cS_1, \cS_3 $ and those subtrees in
the definition of the $N(\G_i^*,\varepsilon)$.
This covers $\Gamma $ and so
\begin{equation}
 \# \cS_2 \le \sum_i N(\Gamma_i^*,\varepsilon) \le 3 \sum_i
M(\Gamma_i^*,\varepsilon) + 3 \# \partial K
\end{equation}
by Lemma 3.9. Let $ \G_i^*(j)$ be the subtrees in the definition
of $M(\Gamma_i^*,\varepsilon)$. Then
\[
\varepsilon \# \cL_2 \le \sum_{i,j} A(\G_i^*(j)) \le
\sum_{i,j}\|u\|_{p',\G_i^*(j)} \|v\|_{p,\G_i^*(j)} \le
\|u\chi_{\Gamma \setminus K}\|_{p'}
\|v\chi_{\Gamma \setminus K}\|_{p} \le \eta^2.
 \]
Since $\# \cL_1\le M(K,\varepsilon)$
\[ M(\G, \varepsilon) \le M(K, \varepsilon) + \# \cL_2 + \# \cL_3,
\]
\[ 0< \varepsilon [M(\G, \varepsilon) - M(K, \varepsilon)] \le
\varepsilon ( \# \cL_2 + \# \cL_3) \le \eta^2 + \varepsilon \# \partial K.
\]
Then, by Lemma 5.1,
\[
0 \le  \limsup_{\varepsilon \to 0} \varepsilon M(\G,\varepsilon) -
 \alpha_p \int_K |uv|  \le
 \eta^2.
 \]
 Now let $K \to \Gamma$ ($\eta \to 0$) to get
 \[ \limsup_{\varepsilon \to 0} \varepsilon M(\G,\varepsilon)
 =
 \alpha_p \int_{\Gamma} |uv|.
 \]

We get the same for $\liminf $ and so
 \[ \lim_{\varepsilon \to 0} \varepsilon M(\G,\varepsilon)
 =
 \alpha_p \int_{\Gamma} |uv|.
 \]
Since $N(K,\varepsilon) + \# \cS_2 + \# \cS_3 \ge N(\G,\varepsilon),$
\[
0\le N(\G,\varepsilon)-N(K,\varepsilon)\le \# \cS_2 +\# \cS_3 \le
3\sum_iM(\G_i^*(j),\varepsilon)+3\# \partial K + \# \partial K
\]
by (5.1). Hence, as before,
\[
\lim_{\varepsilon \to 0} \varepsilon N(\G,\varepsilon)=\alpha_p
\int_{\Gamma}|uv|.
\]
{$\Box$}

\bc Let $1<p<\infty$, $u \in L^{p'}(\G)$ and $v\in
L^p(\G)$. Then
\[
\lim_{n \rightarrow \infty}na_n(T) = \alpha_p \int_{\Gamma}|u||v|.
\]
\ec
{\bf Proof.}
Note that the application of Theorem 3.18 to Theorem 5.3 implies
that $\lim_{n \rightarrow \infty} a_n(T)=0$ and hence that $T$ is
compact.

Let $\{\G_l \}_{l=1}^{\infty}$ be as in the proof of Lemma 5.2, and set
$T_l=T_{a,\G_l}$ for some $a\in \G$. Since $\G_l$ is compact then from
Lemma 3.19 and Theorem 5.1 we have
\[
\lim_{n \to \infty} n a_n(T_l) = {\alpha_p} \, \int_{\G_l} |u||v|.
\]
An operator of rank $<n$ on $\G_l$ can be considered as the
restriction to $ \G_l$ of such an operator on $\G$ and also
$ \|T|L^p(\G)\rightarrow L^p(\G)\| \ge \|T_l|L^p(\G_l)\rightarrow
L^p(\G_l)\| $ if $T=T_a$ and $a \in \G_l$. It follows that $a_n(T)
\ge a_n(T_l)$ and so
\[
\liminf_{n \to \infty} n a_n(T) \ge {\alpha_p} \, \int_{\G_l}
|u||v|.
\]
But, we know from Lemma 5.3 that
\[
  \lim_{\varepsilon \to 0_+}\varepsilon N(\G,\varepsilon)
= \alpha_p \int_{\G} |u| |v|.
\]
and so, by Lemma 3.16

\[\limsup_{n\to \infty} n a_n(T) \le  \, \alpha_p
\int_{\Gamma} |u| |v|.
\]
Hence the corollary is proved.{$\Box$}

{ {\bf 5.2 The cases $p=\infty$ and $p=1$.} \hfill }

The analogies of Lemma 5.1 for $p=\infty$ and $p=1$ are,
respectively

\bl Let $K$ be a compact subtree of $\G,$ and $p=\infty$. Then
\[ {1 \over 2} \int_K |u| |v_s| \le \liminf_{\eps \to 0_+} \eps
N(K,\eps,u,v) \le \limsup_{\eps \to 0_+} \eps N(K,\eps,u,v) \le
{3 \over 2} \int_K |u| |v_s| \]
where $v_s$ is defined in Definition 5.4.
\el

\bl Let $K$ be a compact subtree of $\G$ and $p=1$. Then
\[{1 \over 2} \int_K |u_s| |v| \le \liminf_{\eps \to 0_+} \eps
N(K,\eps,u,v) \le \limsup_{\eps \to 0_+} \eps N(K,\eps,u,v) \le
{3\over 2} \int_K |u_s| |v|. \]
\el

Both lemmas follow from the results for intervals in [6] since $K$
is a finite union of intervals. Lemmas 5.5 and 5.6 yield, as in
Lemma 5.2,

\bl  For $p=\infty$
\[{1\over 2} \int_{\G} |u| |v_s| \le \liminf_{\eps \to 0_+} \eps
N(\G,\eps,u,v) \]
and for $p=1$
\[{1\over 2} \int_{\G} |u_s| |v| \le \liminf_{\eps \to 0_+} \eps
N(\G,\eps,u,v). \]
\el

\bl Let $u\in L^{p'}(\G)$ and $v\in L^p(\G)$. Then for $p=\infty$
\[{1\over 2} \int_{\G} |u| |v_s| \le \liminf_{\eps \to 0_+} \eps
N(\G,\eps,u,v) \le \limsup_{\eps \to 0_+} \eps
N(\G,\eps,u,v) \le {3\over 2} \int_{\G} |u| |v_s| \]
and for $p=1$
\[{1\over 2} \int_{\G} |u_s| |v| \le \liminf_{\eps \to 0_+} \eps
N(\G,\eps,u,v) \le \limsup_{\eps \to 0_+} \eps N(\G,\eps,u,v) \le
{3\over 2} \int_{\G} |u_s| |v|.
\]
\el
{\bf Proof.} Let
$p=\infty$. We need only prove the last inequality. Let
$\{\G_l\}_{l=1}^{\infty}$ be compact subtrees of $\G$ which are
such that
\[ \left|\int_{\G} |u| |v_s|- \int_{\G_l}|u| |v_s| \right| \le
{1 \over l} \]
and
\[ \|u\|_{1,{\G \setminus \G_l}} \le {1\over l}.
\]
Fix $l\in \N$. There exist intervals $W(j)$ in $\G_l$ and step functions
$u_n,v_n$ on $\G_l$,
\[
u_{n}=\sum_{j=1}^m \xi_j \chi_{W(j)}, \qquad
v_{n}=\sum_{j=1}^m \eta_j \chi_{W(j)},
\]
which are such that
\[\|u-u_n \|_{1,\G_l} <{1\over n}, \quad \int_{\G_l} |u(t)|
(v_n(t)-v_s(t)) dt <{1\over n}
\]
and $\|v_s\|_{\infty,\G} \ge
v_n(t) \ge v_s(t)$ on $\G_l$; cf [6, Theorem 4.7].

Let $M:=M(\G,\varepsilon)$ and let $\{\G_i^M\}_{i=1}^M$ be a maximal set
of subtrees of $\G$ in the definition of $M$ (see Definition 3.11).
Then, because $\G_l$ is a compact subtree of $\G$, we have
$ M-2m-\#V(\G_l)-\#\partial\G_l \le \#\K,$ where
\[
\K:=\{\G_j :  \G_j\in\{\G_k^M\}_{k=1}^M, \mbox{ and there exists
} i \mbox{ such that } W(i) \supset \G_j \mbox{ or } \G_j \subset
\G\setminus \G_l\}. \]
On using Lemmas 4.6 and 4.7, we have

\begin{eqnarray*}
\eps(M-2m &-& \#V(\G_l)-\#\partial\G_l)  \le  \sum_{k\in \K}
A(\G_k^M,u,v) \\
& \le & \sum_{k\in \K} \left( A(\G_k^M,u_n,v_n)+ \left[
A(\G_k^M,u,v)-A(\G_k^M,u_n,v) \right] \right.\\
& & \qquad + \left. \left[ A(\G_k^M,u_n,v)-A(\G_k^M,u_n,v_n) \right]
\right)\\
& \le & {1\over 2} \sum_{j=1}^m |\xi_j| |\eta_j| |W(j)| \\
& & +
\sum_{j=1}^M \left( \|u-u_n\|_{1,\G^M_j} \|v\|_{\infty,\G_J^M}
\right)\\
& & + \sum_{k\in \K; \exists i, W(i) \supset \G_k^M}
\left[ A(\G_k^M,u_n,v)-A(\G_k^M,u_n,v_n) \right] \\
& & + \sum_{k\in \K; \G_k^M\subset \G\setminus \G_l}
\left[ A(\G_k^M,u_n,v)-A(\G_k^M,u_n,v_n) \right] \\
& \le &  {1 \over 2} \sum_{i=1}^m |\xi_j| |\eta_j| |W(j)|+
\|u-u_n\|_1 \|v\|_{\infty} \\
& & +\sum_{k\in K;\exists i, W(i) \supset \G_k^M } \left[{\alpha
\over 2} \int_{\G^M_k} (v_n -v_s)|\xi_i| dt + {|\xi_i| |\eta_i|
\over 2\alpha} |\G_k^M| \right] \\
& & +  \sum_{k\in \K, \G_k^M \subset \G\setminus \G_l}
\left[A(\G_k^M,0,v)-A(\G_k^M,0,0) \right] \\
& \le &  {1 \over 2} \sum_{i=1}^m |\xi_j| |\eta_j| |W(j)|+({1\over n} + {1\over
l})\|v_s\|_{\infty} \\
& & + {\alpha
\over 2} \int_{\G} (v_n -v_s)|u_n| dt\\
& &  +{1\over
2\alpha} \int_{\G} |u_n| |v_n| dt \\
& \le & {1\over 2} \int_{\G} |u| |v_s| + c\left(\alpha {1\over
n}+{1\over \alpha} +{1\over n} +{1\over l} \right)
\end{eqnarray*}
for some constant $c$ independent on $\eps$. We therefore
conclude that
\[
{1\over 3} \limsup_{\eps\to 0_+}\eps N(\G,\eps,u,v) \le
{1\over 2} \int_{\G} |u| |v_s|+ K(\alpha {1\over n}+ {1\over
\alpha}),\]
whence
\[{1\over 3} \limsup_{\eps \to 0_+} \eps N(\G,\eps,u,v) \le
{1\over 2} \int_{\G} |u| |v_s|. \]

The case $p=1$ is similar.
{$\Box$}

From Lemmas 5.8 and 3.18 we derive

\bl
Let $u\in L^{p'}(\G)$ and $v\in L^p(\G)$. Then
for $p=\infty$
\[
{1 \over 6} \int_{\G} |u| |v_s| \le \liminf_{n \to \infty} n a_n(T)
 \le \limsup_{n \to \infty } n a_n(T)  \le 3  \int_{\G} |u| |v_s|
\]
and for $p=1$
\[
{1\over 6} \int_{\G} |u_s| |v| \le \liminf_{n \to \infty} n a_n(T)
 \le \limsup_{n \to \infty } n a_n(T)  \le 3  \int_{\G} |u_s| |v|.
\]
\el

\section{The main results for $1<p<\infty$.}

We suppose throughout this section that $T:=T_a$ and $\G:=\G_a$
for some $a\in \G$. Also we write $\|T_K\|$ for $\|T|L^p(K)
\rightarrow L^p(K)\|$, for any $K \subseteq \G$.

With $U(x):= \int_a^x |u(t)|^{p'}dt$ $(x\in \G)$
we define $Z_k$ to be the closure of
\begin{equation}
\{x: x\in \G, 2^{kp'\over p} \le U(x) <
2^{(k+1)p'\over p} \}. \label{6.1}
\end{equation}
Here $k$ may be any integer if $u\in L^{p'}_{loc}(\G)\setminus
L^{p'}(\G)$, while, if $u\in L^{p'}(\G)$, $2^k\le \|u\|^p_{p',\G}$; we
refer to these values of $k$ as the {\em admissible} values.

We have that $Z_k=\bigcup_{i=1}^{n_k} Z_{k,i}$, where the
$Z_{k,i}$ are the connected components of $Z_k$. Corresponding to
each admissible $k$ we set
\begin{equation}
\sigma_{k,i}^p:=2^k \mu(Z_{k,i})
 \mbox{ for } i\in \{1, \dots ,n_k\} \label{7.3}
\end{equation}
and
\begin{equation}
\sigma_{k}^p:= 2^k \mu(Z_{k}).
\label{7.4}
\end{equation}
For non--admissible $k$ we set $\sigma_k=0$. We also set
$\sigma_{k,i} =0 $ for $i\notin \{1,\dots, n_k \} $.

Let
\beq B_{k,i}:= \# \partial Z_{k,i} -1 \label{6.4};
\eeq
that is, $B_{k,i} $ is the number of boundary points of $Z_{k,i}$ excluding its root.

\bl
\beq \sup_{k\in \Z} \max_{1\le j \le n_{k}} \sigma_{k,i} \le \|T\|. \label{6.5}
\eeq
\el
{\bf Proof.} This follows from [7, Proposition 5.1], which asserts that
\[ \sup_{x\in \Gamma} \|u\chi_{(a,x)}\|_{p'} \|v \chi_{(a,x)^c} \|_p \le \|T\|, \]
where $(a,x)^c=\{y \in \Gamma: x \preceq_a y\}$. For then, by
(\ref{6.1}),
\begin{eqnarray*}
\|T\| &\ge& \sup_{x\in \G} U(x)^{1/p'}\left(\int_{y\succeq
x}|v(y)|^pdy\right)^{1/p}\\
&\ge& \sup_{k,i} 2^{k/p} \mu(Z_{k,i})^{1/p}\\
 &=& \sup_{k,i} \sigma_{k,i} .
\end{eqnarray*}
{$\Box$}

\bl
 Let $\Gamma'$ be a subtree of $\Gamma =\Gamma_a$
and $b=b(\Gamma')$ the nearest point of $\Gamma'$ to $a$. Then, for any $c>4$,
there exist
$X=X(\Gamma') \in \I_b(\Gamma')$ and $ k'=k'(\Gamma') \in \Z$ such that, with
$T'=T_{\Gamma'}$ and $Y=Y(\Gamma')=\Gamma' \setminus X,$
\beq \|T'\| \le 2^{2/p} c \{\sum_{i \in S} 2^{k'} \mu(Y \cap Z_{k',i}) \}^{1/p}, \label{6.6}
\eeq
where $S=S(\Gamma')=\{i:\mu(Y\cap Z_{k',i})>0\}.$
\el
{\bf Proof.} From Theorem 2.4, for $c>4$, there exists $X\in \I_b(\Gamma')$ such that
\begin{eqnarray*}
\|T'\| & \le & c {\mu(Y)^{1/p} \over \alpha_X} \\
& \le & c
\min_{t \in \partial X\setminus\{b\}} \left( \int^t_b |u|^{p'} dx \right)^{1/p'} \mu(Y)^{1/p} \\
& \le & c \min_{t \in \partial X\setminus\{b\}} [U(t) - U(b)]^{1/p'}
[\sum_{i,k} \mu(Y \cap Z_{k,i})]^{1/p} \\
& \le & c \min_{t \in (\partial X\setminus\{b\})\cap Z_{\gamma_o}}
[U(t) - U(b)]^{1/p'}
[\sum_{i,k} \mu(Y \cap Z_{k,i})]^{1/p},
\end{eqnarray*}
where $\gamma_0=\min \{k: \mu(Y\cap Z_k)>0\}$. Since
$\mu(\G)<\infty$, we may assume that $Y$ is compact and hence
\[\max_{k\ge \gamma_0} \sum_{i=1}^{n_k} 2^k \mu(Y\cap Z_{k,i}) \]
is attained, and so
\begin{eqnarray*}
\sum_{k,i} \mu(Y\cap Z_{k,i}) & = & \sum_{k \ge \gamma_0} 2^{-k} \sum_{i=1}^{n_k}
2^k \mu(Y\cap Z_{k,i}) \\
& \le & 2^{1-\gamma_0} \sum_{i=1}^{n_{k'}} 2^{k'} \mu(Y\cap Z_{k',i})
\end{eqnarray*}
say, for some $k' \ge \gamma_0.$ Hence
\[ \|T'\| \le c[ 2^{(\gamma_0+1)p'/p}-2^{\gamma_0 p'/p} ]^{1/p'}
[2^{1-\gamma_0} \sum_{i=1}^{n_{k'}} 2^{k'} \mu(Y\cap Z_{k',i}) ]^{1/p}, \]
whence (\ref{6.6}).
{$\Box$}

\bl Let $\{ \Gamma_i \}_{\Alpha}$ be a finite set of non-overlapping subtrees of $\Gamma$ and set
$T_l=T_{\Gamma_l}$. Then,
\beq \sum_{l \in \Alpha} \|T_l\|^q \le (2^{2/p+2})^q \sum_{(k,i)\in \eta} B_{k,i}^{q/p'} \sigma_{k,i}^q
\qquad \mbox{ if } 1\le q \le p  \label{6.7}
\eeq
and
\beq \sum_{l \in \Alpha} \|T_l\|^q \le (2^{2/p+2})^q \sum_{k\in \eta_0}  \sigma_{k}^q \qquad
\mbox{ if } p \le q < \infty , \label{6.8}
\eeq
where $\eta $ and $\eta_0$ are finite sets.
\el
{\bf Proof.} Let $\Gamma_{\lambda} \in \{\Gamma_l\}_{\Alpha}$, and, in the notation of Lemma 6.2,
set $b_l=b(\Gamma_l),$ $k_l=k'(\Gamma_l)$, $Y_l=Y(\Gamma_l)$ and $S_l=S(\Gamma_l)$. There are two cases to
consider for $\Gamma_{\lambda}$:

\begin{description}
\item[(i)] $b_{\lambda} \in Z_{k_{\lambda}}.$ In this case $b_{\lambda} \in Z_{k_{\lambda},i_{\lambda}},$
$S_{\lambda}=S(\Gamma_{\lambda})=\{ i_{\lambda} \}$ and, for any
$c>4$,
\[ \|T_{\lambda}\| \le (2^{2/p} c) \sigma_{k_{\lambda},i_{\lambda}}. \]
\item[(ii)] $b_{\lambda} \not \in Z_{k_{\lambda}}.$ Denote by $\Lambda$ the subset of $\Alpha$ which is such
that for $l\in \Lambda$, $b_l \in Z_{k_{\lambda},i_l}$ for some unique $i_l \in S_{\lambda}$ and
so $S_l=\{i_l\}$.
\end{description}

Set $\Lambda_i=\{l\in \Lambda: i_l=i\}.$ Then, by (\ref{6.6}), for $q\ge 1$,

\begin{eqnarray}
\sum_{l\in \Lambda} \|T_l\|^q & = & \sum_{i\in S_{\lambda}} \sum_{l \in \Lambda_i} \|T_l\|^q \nonumber \\
& \le & (2^{2/p}c)^{q} \sum_{i\in S_{\lambda}} \sum_{l \in \Lambda_i}
\left\{ 2^{k_{\lambda}} \mu(Y_\l \cap Z_{k_{\lambda},i}) \right\}^{q/p}  \label{6.9} \\
& \le & (2^{2/p}c)^{q} \sum_{i\in S_{\lambda}} \left\{ \sum_{l \in \Lambda_i}
\left[ 2^{k_{\lambda}} \mu(Y_\l \cap Z_{k_{\lambda},i}) \right]^{1/p} \right\}^{q}  \nonumber \\
& \le & (2^{2/p}c)^{q} \sum_{i\in S_{\lambda}} \left[ \left\{ \sum_{l \in \Lambda_i}
 2^{k_{\lambda}} \mu(Y_\l \cap Z_{k_{\lambda},i} ) \right\}^{1/p}
 \left\{ \sum_{l\in \Lambda_i} 1 \right\}^{1/p'} \right]^q  \nonumber \\
& \le & (2^{2/p}c)^{q} \sum_{i\in S_{\lambda}} \left( B_{k_{\lambda},i}^{1/p'}
\sigma_{k_{\lambda},i}\right)^q. \label{6.10}
\end{eqnarray}
Also in case $(ii)$, from (\ref{6.6}),
\[\|T_{\lambda} \| \le (2^{2/p}c)(\sum_{i\in S_{\lambda}} \sigma_{k_{\lambda},i}^p)^{1/p}.\]
Hence, if $1\le q \le p$,
\begin{eqnarray}
\|T_{\lambda} \|^q & \le & (2^{2/p}c)^q (\sum_{i\in S_{\lambda}} \sigma_{k_{\lambda},i}^q) \nonumber \\
& \le & (2^{2/p}c)^q \sum_{i\in S_{\lambda}} \left(
 B_{k_{\lambda},i}^{1/p'}
\sigma_{k_{\lambda},i}\right)^q. \label{6.11}
\end{eqnarray}

If $q\ge p$, then from (\ref{6.9}),
\begin{eqnarray}
\sum_{l\in \Lambda} \|T_{l} \|^q & \le & (2^{2/p}c)^q \sum_{i\in S_{\lambda}}
\left\{ \sum_{l\in \Lambda_i} 2^{k_{\lambda}} \mu(Y_l\cap Z_{k_{\lambda},i}) \right\}^{q/p} \nonumber \\
& \le & (2^{2/p}c)^q \sum_{i\in S_{\lambda}} \sigma_{k_{\lambda},i}^q \nonumber \\
& \le & (2^{2/p}c)^q \left( \sum_{i\in S_{\lambda}} \sigma_{k_{\lambda},i}^p \right)^{q/p} \nonumber \\
& \le & (2^{2/p}c)^q \sigma_{k_{\lambda}}^q. \label{6.12}
\end{eqnarray}

Also, by (\ref{6.6}),
\begin{eqnarray}
\|T_{\lambda} \|^q & \le &
(2^{2/p}c)^q \left\{ \sum_{l\in S_{\lambda}} 2^{k_{\lambda}}
\mu(Y_{\lambda}\cap Z_{k_{\lambda},i}) \right\}^{q/p} \nonumber \\
& \le & (2^{2/p}c)^q \sigma_{k_{\lambda}}^q. \label{6.13}
\end{eqnarray}
The lemma follows from (\ref{6.10})-(\ref{6.13}) since $c>4$ is arbitrary.
{$\Box$}

\bt For $1<p<\infty$, let $u,v$ satisfy (\ref{2.1}) and suppose that
$ B_{k,i}^{1/p'} \sigma_{k,i} \in l^1(\Z \times \N) $. Then
\beq
\lim_{\varepsilon \to 0} \varepsilon M(\G,\varepsilon)=
 \alpha_p \int_{\G} |u| |v|, \label{6.14}
\eeq
\beq \lim_{\varepsilon \to 0} \varepsilon N(\G,\varepsilon,)=
 \alpha_p \int_{\G} |u| |v|, \label{6.15}
\eeq
and
\beq \lim_{n\to \infty} n a_n(T) = \alpha_p \int_{\G} |u| |v|. \label{6.16}
\eeq

\et
{\bf Proof.}  Given $\eta > 0$, we choose $l$ to be such that
$ \sum_{k\ge l} \sum_{i=1}^{n_k} B_{k,i}^{1/p'} \sigma_{k,i} \le \eta $ and set
$K= \bigcup_{k\le l} Z_k $. Then, in the notation of the proof
of Theorem 5.3, we have by Corollary 3.3, that
\begin{eqnarray*}
\varepsilon \# \Alpha_2 & \le &
\sum_{i,j} A(\G_i^{*(j)}) \le \sum_{i,j} \|T_{\G_i^*(j)} \| \\
 & \le & c \sum_{k\ge l} \sum_{i=1}^{n_k} B_{k,i}^{1/p'} \sigma_{k,i}
\end{eqnarray*}
for some positive constant $c$, by Lemma 6.3 with $q=1$.
The proofs of the first two
identities then follow that of Theorem 5.3. Theorem 3.18 and Lemma 3.16
complete the proof.

Note that the convergence of
$\sum_{k\in \Z} \sum_{i=1}^{n_k} B_{k,i}^{p'} \sigma_{k,i} < \infty$
implies that $T$ is compact. To see this, let $K_j=\cup_{k\le j} Z_k$, and so
\begin{eqnarray*}
(T_{K_j}f)(x) & = & v(x) \chi_{K_j}(x) \int_a^x f(t) u(t) \chi_{K_j}(t) dt \\
& = & v(x) \chi_{K_j}(x) \int_a^x f(t) u(t) dt.
\end{eqnarray*}
Then
\[(T-T_{K_j})f(x)= v(x) \chi_{\Gamma \setminus K_j}(x) \int_a^x f(t) u(t) dt \]
and, by Lemma 6.2, for some $k'>j$
\begin{eqnarray*}
\|T - T_{K_j} \|  & \le & 2^{2/p} c \left\{ \sum_{i\in S} 2^{k'} \mu(Z_{k',i})\right\}^{1/p} \\
& \le & 2^{2/p}c \sum_{i\in S} \sigma_{k',i} \le 2^{2/p} c
\sum_{i\in S} B_{k',i}^{1/p'} \sigma_{k',i}.
\end{eqnarray*}
Thus $\|T-T_{K_j}\|\to 0 $ as $j \to \infty $, and $T$ is compact since the
$T_{K_j}$ are compact.{$\Box$}

\bt Let $ 1<q \le p$. Then,
for some positive constant $c$,
\beq \| \{ a_n(T)\} \|_{l^q(\N)} \le c \| B_{k,i}^{1/p'} \sigma_{k,i} \|_{l^q(\Z \times \N) }
\mbox{ if } 1<q \le p \label{6.17}
\eeq
\et
{\bf Proof:} Let $\Gamma_l, l=1, 2, \dots, M(\Gamma,\varepsilon)$ be a maximal set of subtrees
of $\Gamma$ from the definition of $M(\Gamma, \varepsilon)$, so that $A(\Gamma_l)>\varepsilon$.
Then, from (\ref{6.7}) and Corollary 3.3, for $1\le q \le p$,
\[ \varepsilon^q M(\Gamma,\varepsilon) \le \sum_l \|T_l\|^q \le
c\|B_{k,i}^{1/p'} \sigma_{k,i} \|^q_{l^q(\Z \times \N)}.
\]
Since $a_{3M(\Gamma, \varepsilon)+4}(T) \le 2 \varepsilon $ by Lemma 3.9 and Theorem 3.18, it
follows that
\begin{eqnarray*}
\# \{m: a_m(T)>t \} & \le & 3M(\Gamma, t/2) +4 \\
& \le & c M(\Gamma, t/2)
\end{eqnarray*}
for $c\ge 7$. Thus
\beq
\#\{m:a_m(T)>t \} \le c t^{-q} \|B^{1/p'}_{k,i} \sigma_{k,i}
\|^q_{l^q(\Z \times \N)}. \label{6.18}
\eeq
We now proceed as in the proof of the Marcinkiewicz
Interpolation Theorem (see [11]); we give the proof for
completeness.

Define

\[
v_1:= \left\{
\begin{array}{ll}
v & \mbox{ on } Z_{k,i} \mbox{ if } B^{1/p'}_{k,i} \sigma_{k,i} \le t/2, \\
0 & \mbox{ otherwise,}
\end{array}         \right.
\]
and set $v_2=v-v_1.$ Denote $T, \sigma_{k,i},$ etc by $T(v),
\sigma_{k,i}(v)$ to indicate the dependence on $v$, and set
$T^j=T(v_j), \  j=1,2.$ Then, by [3, Proposition II.2.2],
\[a_{2n-1}(T) \le a_n(T_1) + a_n(T_2) \]
and  so
\[
\{n: a_{2n-1}(T) >t\} \subseteq \{n: a_n(T_1)> t/2\} \cup \{n:a_n(T_2)
>t/2\} \]
and
\beq
\#\{n: a_{2n-1}(T) >t\} \le \#\{n: a_n(T_1)> t/2\} + \#\{n:a_n(T_2)
>t/2\}.  \label{6.19}
\eeq

Set $S_{k,i}=B^{1/p'}_{k,i} \sigma_{k,i}$, and let $1<q<q_1$.
Then, on using (\ref{6.18}) and (\ref{6.19}),
\begin{eqnarray*}
\|\{a_{2n-1}(T)\}\|_{l^q(\N)} & = & q\int_0^{\infty} t^{q-1}
\#\{n: a_{2n-1}(T)>t\} dt \\
& \le & cq \int_0^{\infty} t^{q-1} \left\{ t^{-q_1} \sum_{S_{k,i}
\le t/2} S_{k,i}^{q_1} + t^{-1} \sum_{S_{k,i}>t/2} S_{k,i}
\right\} dt \\
& \le & c \sum_{k,i} S_{k,i}^q,
\end{eqnarray*}
whence (\ref{6.17}), since $a_n(T)$ decreases with $n$.
{$\Box $}

\bt
Let $q\in (p,\infty)$. Then, for some positive constant $c$,
\beq
\| \{a_n(T)\} \|_{l^q(\N)} \le c \| \{\sigma_k\} \|_{l^q(\Z)}.
\label{6.20}
\eeq
\et
{\bf Proof.} Let $\{\Gamma_l \}_1^{M(\Gamma, \varepsilon)}$ be as
in the proof of Theorem 6.5 and define
\[
F_j=\{\Gamma_l: k'(\Gamma_l)=  k_j \}
\]
in the notation of Lemma 6.2. Then, from the proof of Lemma 6.3,
for $q\in (p,\infty),$
\begin{eqnarray*}
\varepsilon^p \#F_j & \le & \sum_{\Gamma_l \in F_j} \|T_l\|^p \\
& \le & c^p \sigma^p_{k_j}.
\end{eqnarray*}
Thus, with $m_j = [c^p \sigma^p_{k_j}/\varepsilon^p],$
\[M(\Gamma, \varepsilon) = \sum_{j=1}^{\infty} \sum_{m=1}^{m_j}
1 \le \sum_{m=1}^{\infty} \# \{j: c\sigma_{k_j} > m^{1/p}
\varepsilon \}. \]
Hence
\begin{eqnarray*}
\|\{ a_n(T) \}\|_{l^q(\N)} & = & q \int_0^{\infty} t^{q-1}
\# \{n: a_n(T)>t\} dt \\
& \le & c \int_0^{\infty} t^{q-1} M(\Gamma, t/2) dt \\
& \le & c \int_0^{\infty} \sum_{m=1}^{\infty} t^{q-1}
\# \{j: \sigma_j>m^{1/p} t \} dt \\
& \le &  c \int_0^{\infty} \sum_{m=1}^{\infty} m^{-q/p}
\# \{j: \sigma_j>t \} t^{q-1} dt \\
& \le & c \|\{\sigma_k \} \|_{l^q(\Z)}.
\end{eqnarray*}
{$\Box$}

In the next theorem $l^q_{\omega}$ denotes weak-$l^q$, that is,
the space of sequences $\{x_k\}$ such that
\[ \|\{x_k\}\|_{l^q_{\omega}}:= \sup_{t>0} \{t(\#
\{k:|x_k|>t\})^{1/q}\}<\infty. \]

\bt
For some positive constant $c$,

\beq
\|\{ a_n(T) \}\|_{l^q_{\omega}(\N)} \le c \|\sum_{i=1}^{n_k} B_{k,i}^{1/p'} \sigma_{k,i} \|_{l^q_{\omega}(\Z)}
              \mbox{ if } 1<q\le p, \label{6.21}
\eeq
and
\beq
\|\{ a_n(T) \}\|_{l^q_{\omega}(\N)}
\le c \|\{\sigma_{k} \}\|_{l^q_{\omega}(\Z) }
              \mbox{ if } p<q<\infty. \label{6.22}
\eeq

\et
{\bf Proof.} Let $\{\Gamma_l\}_1^{M(\Gamma,\varepsilon)}, F_j$ be as in the proof of Theorem 6.6. Then,
from the proof of Lemma 6.3,
\begin{eqnarray*}
\varepsilon \#F_j & \le & \sum_{\Gamma_l \in F_j} \|T_l\| \\
& \le & c \sum_{i=1}^{n_{k_j}} B^{1/p'}_{k_j,i} \sigma_{k_j,i} =:N_j
\end{eqnarray*}
say. Thus
\begin{eqnarray*}
\varepsilon^q \# \{n:a_n(T)>\varepsilon \} & \le & c\varepsilon^q M(\Gamma, \varepsilon/2) \\
& \le & c \varepsilon^q \sum_{j=1}^{\infty} \sum_{m=1}^{[N_j/\varepsilon]} 1 \\
& \le & c \varepsilon^q \sum_{m=1}^{\infty} \# \{j: N_j>m\varepsilon \} \\
& \le & c \ \sum_{m=1}^{\infty} m^{-q} t^q\# \{j: N_j> t \}
\end{eqnarray*}
and hence (\ref{6.21}). The proof of (\ref{6.22}) is similar, starting from
\[
\varepsilon ^p \#F_j \le c \sigma^p_{k_j}.
\]{$\Box$}

Let us now suppose that the tree $\G$ satisfies the following condition:
\beq B_{k,i} < B < \infty \quad \mbox{for each admissible } k \mbox{ and }
i.
\label{6.23}
\eeq

Then with this condition we can get lower estimates in Theorems 6.5 and 6.6.

We need the following result which is similar to [2,Lemma 20].

\bl Suppose that (6.10) is satisfied. Let $S(\varepsilon) := \{(k,i):
\sigma_{k,i}>\varepsilon \}$. If
$M+1/2 \le \# S(\varepsilon)/4B $, then $a_M(T)>c\varepsilon $,
where $c$ is an absolute constant.
\el
{\bf Proof.}
It is sufficient to prove the result for $S(\varepsilon)$ finite, for
this will imply the result when $ \#S(\varepsilon)= \infty $. The
elements of $S(\varepsilon) $ fall into two subsets according as
$k$ is odd or even. At least one of them, say $S_1(\varepsilon)$,
has cardinal at least half that of $S(\varepsilon)$. Thus, we may
suppose that $\#S_1(\varepsilon) > B$.

Denote by $ \zeta_{k,j}$ the point of $Z_{k,j}$ nearest to $a$,
and define $n(x) := \#\{(k,j):\zeta_{k,j}\succ_a x, (k,j)\in
S_1(\varepsilon)\}$. Let $l$ be a path in $\G$ starting at $a$ and
consisting of edges $(x,y)$ of $\G$, $(x,y)$ at each stage chosen so that
$n(y)$ is as large as possible. Terminate the path at the point $ x=
\zeta_{r,s}$ at which $n(x)=0$. Define $\xi \in l$ by
\[
\xi := \inf \{x \in l:n(x)=n(\zeta_{r-1,j}), \zeta_{r-1,j} \in l\},
\]
( the infimum, which is being taken with respect to the total ordering on
$l$ induced by $\preceq_a $, exists since $n(a)>B$ and
$ n(\zeta_{r-1,j}) \le B$ ). There are two possibilities : (i) $\xi $
may be a point $\zeta_{k,j}$, in which case define $\G_1 := \{x:
x\succ_a \xi \}$, or (ii) $\xi $ may be a vertex of $\G $ joined by a path
$l_1$ to a point $\zeta_{m,n} \succ_a \xi $, where $(m,n) \in
S_1(\varepsilon)$. In the latter case we define $\G_1
:=\{x:x\succ_ay, y\in l\cup l_1\}$. Then, in both cases, the
closure of $\G_1$ is a subtree and so is its complement. Moreover,
$A(\overline{\G_1})>c\varepsilon $, where $c$ is an absolute
constant. For, in case (i), if $b$ is a point of $l$ with
$U(b)=2^{(r-1/2)(p'/p)}$ and $T^1 $ is the restriction of $T$ to
$\G_1$, then, in the notation of the discussion preceding Lemma 3.5,
\[
\|T^1_{b,1}\|, \|T^1_{b,2}\| \ge
(2^{r(p'/p)}-2^{(r-1/2)(p'/p)})^{1/p'}(2^{-r/p}\varepsilon );
\]
this follows from [7, Proposition 5.1] where it is shown that
\[
\|T_a\| \ge \sup_{x\in \G}\|u\chi_(a,x)\|_{p'}\|v\chi_{(a,x)^c}\|_p.
\]
In case (ii) a similar result holds if $b$ is a point of
$l\cup l_1 $ with $U(b)=2^{(t-1/2)(p'/p)}$ and $t$ the greater of
$r,m$. The lower bound for $A(\overline{\G_1})$ is then a
consequence of Lemma 3.5. Note also that $\G_1$ contains at most
$2B$ elements of $S_1(\varepsilon) $.

The result now follows by induction on $\#S(\varepsilon)$ and
Lemma 3.17.
{$\Box $}
\bl
Suppose that (6.10) is satisfied. Then, for all $t>0$,
\[
\# \{(k,i): \sigma_{k,i}>t \} \le 4B \# \{k\in \N: a_k(T)>
ct \} +6B.
\]
\el
{\bf Proof.} From Lemma 6.8,
\begin{eqnarray*}
\# \{k\in \N: a_k(T)> ct \} &\ge & \big[\#S(t)/4B-1/2 \big] \\
& \ge & \#S(t)/4B -3/2,
\end{eqnarray*}
whence the result.
{$\Box $}

\bl
Suppose that (6.10) is satisfied. Then, for all $q>0$,
\[ \| \{ \sigma_{k,i}\} \|^q_{l^q({\Z \times \N} )} \le c_1 \|
\{a_k(T)\} \|^q_{l^q(\N)} + c_2 \| \{ \sigma_{k,i} \}
\|^q_{l^{\infty}({\Z \times \N})} \]
\el
{\bf Proof.}  Let $\lambda = \| \{\sigma_{k,i}\}
\|_{l^{\infty}(\Z \times \N)}$. Then, by Lemma 6.9,
\begin{eqnarray*}
 \|\{\sigma_{k,i}\}\|^q_{l^q(\Z \times \N)} & \le & q
\int_0^{\lambda} t^{q-1} \#\{(k,i)\in \Z \times \N: \sigma_{k,i} >t
\} dt \\
& \le & 4B  q
\int_0^{\lambda} t^{q-1} \#\{ k \in \N; a_k(T) > c \varepsilon
\} dt + 6B \lambda^q \\
& \le & c_1 \|\{ a_k(T) \} \|^q_{l^q(\N)} + c_2 \lambda^q.
\end{eqnarray*}
{$\Box$}

\bt
Let $1<p< \infty $ and suppose that (6.10) is satisfied. Then, for any $q>0$,
there exists a constant $c>0$ such that
\[\|\{ \sigma_{k,i} \}\|_{l^q(\Z \times \N)} \le c \| \{ a_k(T)\}
\|_{l^q(\N)}. \]
\et
{\bf Proof.} By Lemma 6.1,
\[ \|\{\sigma_{k,i} \} \|_{l^{\infty}(\Z \times \N)} \le  \|T\|
=  a_1(T) \le  \| \{a_k(T) \} \|_{l^q(\N)}. \]
The result then follows from Lemma 6.10.
 {$\Box$}

\br
{\bf (i)} It follows from Theorem 6.5 and 6.11 that if (\ref{6.23})is satisfied
and $1<q\le p$, then
\[ \|\{a_n(T)\} \|_{l^q(\N)} \asymp \|\{ \sigma_{k,i} \} \|_{l^q(\Z \times \N)}.\]
For $q>p$, we have from Theorems 6.11 and 6.6
\begin{eqnarray}
c_1 \|\sigma_{k_i}\|_{l^q(\Z \times \N)} & \le & \| \{a_n(T) \}\|_{l^q(\N)} \le c_2 \|\sigma_k\|_{l^q(\Z)}
\nonumber \\
& \le & \|\sigma_{k,i} \|_{l^p(\Z \times \N)}. \label{6.24}
\end{eqnarray}

{\bf (ii)} Naimark and Solomyak [8] take $u=1$, and in [8,(4.8)] they make
the assumption that, for every edge $\langle y,z \rangle \in
E(\Gamma),$
\beq
\mu_1 \le |z|/|y| \le \mu_2, \qquad 1<\mu_1 \le \mu_2, \label {6.25}
\eeq
where $|y|,|z|$ denote the lengths of the paths from the root of $\G$
to $y,z $ respectively.
Let $y_j\in V(t),$ $j=0,1, \dots ,$ and suppose that $|y_0|\le 2^k$ and $|y_1| \ge 2^k$. Then (\ref{6.25})
implies that
\[ |y_n|\ge \mu_1^{n-1} |y_1| \ge \mu_1^{n-1}2^k \ge 2^{k+1} \]
if $n \ge 1+ \log 2 / \log \mu_1$. Hence, if each vertex has constant branching number $b$ (ie. degree $b+1$),
then
\[ B_{k,i} \le b^{[\log 2 / log \mu_1 +1]} \]
and hence (\ref{6.23}) is satisfied.

{\bf (iii)} Theorem 4.1 in [9] is valid under assumptions made on a
sequence $\{\eta_j\}$ which is defined as follows : for any
partition $ \Xi $ of $\G$ into a countable union of non-overlapping
segments $I_j = \langle y_j,z_j\rangle $,
\[
\eta_j :=|z_j|\int_{I_j}v^2 dt.
\]
Note that in our notation, the case $p=2, u=1 $ is what
is considered in [9]. It is proved in [9, Theorem 4.1 (i)] that (6.15)
for $p=2$ holds if, for some $\Xi$, $\{\eta_j\}\in l_{1/2}. $

Choose $\Xi = \cup _{k\in \N}\Xi_k $, where $ \Xi_k $ is a
partition of $Z_k$. Then,
\begin{eqnarray*}
\sigma_{k,j}^2 &=& 2^k\sum_{I_s\subset Z_{k,j}} \int_{I_s}v^2 dt \le \sum_{I_s\subset Z_{k,j}} \eta_s \\
&\le & \left( \sum_{I_s\subset Z_{k,j}} \eta_s^{1/2} \right)^2
\end{eqnarray*}
and
\[
\sum_{k,j} \sigma_{k,j} \le \sum_s \eta_s^{1/2}. \]
Thus, if (\ref{6.23}) is satisfied,
\beq
\sum_{k,j} B^{1/2}_{k,j} \sigma_{k,j} \le  B^{1/2} \sum_s \eta_s^{1/2}. \label{6.26}
\eeq

In the reverse directions we have
\[\sigma_{k,j}^2 \ge 1/2 \sum_{I_s \subset Z_{k,j}} \eta_s \]
and so
\[
 \sigma_{k,j} \ge 2^{-1/2} (\sum_{I_s \subset Z_{k,j}} \eta_s^{1/2})
 (\sum_{I_s \subset Z_{k,j}} 1 )^{-1/2}.
\]
Therefore
\[ B_{k,j}^{1/2} \sigma_{k,j} \ge  2^{-1/2} \left( {B_{k,j} \over \sum_{I_s \subset Z_{k,j}} 1}
\right)^{1/2} \sum_{I_s \subset Z_{k,j}} \eta_s^{1/2}.\]
If
\beq \inf_{k,j} \left( {B_{k,j} \over \sum_{I_s \subset Z_{k,j}} 1}
\right) =: c >0 \label{6.27}
\eeq
then
\beq
\sum_{k,j} B_{k,j}^{1/2} \sigma_{k,j} \ge (c/2)^{1/2} \sum_{s}
\eta_s^{1/2}.
\label{6.28}
\eeq

The condition (\ref{6.27}) is satisfied if the tree $\Gamma$ is, in the terminology of [9],
$b$ regular of type $(b,2)$ and $\Xi$ consists of edges of $\Gamma$. This means that every vertex of
$\Gamma$ has fixed branching number $b$, and any edge $\langle y,z\rangle$ of the k-th generation
is such that $|y|=2^k, |z|=2^{k+1}.$ Hence, in this case, (\ref{6.27}) is
satisfied with $c=1$.
\er

{\bf Acknowledgment.}
J. Lang wishes to record his gratitude to the Royal Society and
NATO for support to visit the School of Mathematics at Cardiff
during 1997/8, under their Postdoctoral Fellowship programmes. He
also thanks the Grant Agency  of the Czech Republic for partial support
under grants  \ 201/96/0431 and  \ 201/98/P017.

\vfil\eject\noindent

\begin{center}
{\bf REFERENCES}
\end{center}
\begin{enumerate}
\item
D.E.Edmunds, W.D.Evans and D.J.Harris. Approximation numbers of
certain Volterra integral operators. {\em J. London Math. Soc.}
(2) 37 (1988), 471--489.
\item
D.E.Edmunds, W.D.Evans and D.J.Harris. Two--sided estimates of
the approximation numbers of certain Volterra integral
operators. {\em Studia Math.} 124 (1) (1997), 59--80.
\item
D.E.Edmunds and W.D.Evans, Spectral Theory and Differential
Operators, {\em Oxford Univ. Press, Oxford,} 1987.
\item
D.E.Edmunds, P.Gurka and L.Pick. Compactness of Hardy--type
integral operators in weighted Banach function spaces. {\em
Studia Math.} 109 (1) (1994), 73--90.
\item
W.D.Evans and D.J.Harris, Fractals, trees and the Neumann Laplacian.
{\em Math. Ann.} 296 (1993), 493-527.
\item
W.D.Evans, D.J.Harris and J.Lang, Two--sided estimates for the
approximation numbers of Hardy--type operators in $L^{\infty}$
and $L^1$. {\em Studia Math.} 130 (2) (1998), 171-192.
\item
Evans, W.D.; Harris, D.J.; Pick, L.
Weighted Hardy and Poincaré inequalities on trees. {\em J. Lond. Math. Soc.}
52 (2)(1995), 121-136.
\item
K.Naimark and M.Solomyak, Eigenvalue Estimates for the weighted Laplacian on
metric trees. Preprint.
\item
J.Newman and M.Solomyak, Two--sided estimates of singular values
for a class of integral operators on the semi--axis, {\em
Integral Equations Operator Theory}, 20 (1994), 335--349.
\item
B.Opic and A.Kufner, Hardy--type Inequalities, {\em Pitman Res.
Notes Math. Ser. 219, Longman Sci. \& Tech., Harlow,} 1990.
\item
E.M.Stein, Singular Integrals and Differentiability Properties of Functions,
{\em Princeton Univ. Press, Princeton,} 1970
\end{enumerate}
\vspace{0.50in}
W.~D.~Evans, D.~J.~Harris \\
School of Mathematics \\
Cardiff University \\
Senghennydd Road \\
Cardiff CF24 4YH \\
Wales, UK \\
E-mail : EvansWD@cardiff.ac.uk\\
 \\
J.~Lang \\
Mathematics Department\\
202 Mathematical Sciences Bldg\\
University of Missouri\\
Columbia, MO 65211 USA\\
(permanent address: Math. Inst. AV CR, Zitna 25, Prague 1, Czech Republic)\\
E-mail : langjan@math.missouri.edu

\end{document}